\documentclass[leqno,12pt]{article}
\usepackage[utf8]{inputenc}
\usepackage[T1]{fontenc}
\usepackage{amsmath}
\usepackage{amsfonts}
\usepackage{amsthm}
\usepackage{color}

\newcommand{\red}{\color{red}\tt}
\newcommand{\blue}{\color{blue}}

\newlength{\oldparindent}
\newcommand{\cL}{{\mathbb {L}}}

\newcommand{\bpf}{\begin{preuve}}
\newcommand{\epf}{ \end{preuve} \medskip}

\newcommand{\benum}{\begin{enumerate}}
\newcommand{\eenum}{\end{enumerate}}

\newcommand{\bitem}{\begin{itemize}}
\newcommand{\eitem}{\end{itemize}}

\newcommand{\brmq}{\begin{rmq}}
\newcommand{\ermq}{\end{rmq}}

\newcommand{\brmqs}{\begin{rmqs}}
\newcommand{\ermqs}{\end{rmqs}}

\newcommand{\bapp}{\begin{application}}
\newcommand{\eapp}{\end{application}}

\newcommand{\bapps}{\begin{applications}}
\newcommand{\eapps}{\end{applications}}

\newcommand{\bdefi}{\begin{definition}}
\newcommand{\edefi}{\end{definition}}

\newcommand{\beq}{\begin{equation}}
\newcommand{\eeq}{\end{equation}}

\def\bpm{\begin{pmatrix}}
\def\epm{\end{pmatrix}}

\newcommand{\bcas}{\begin{cases}}
\newcommand{\ecas}{\end{cases}}

\newcommand{\bex}{\begin{exemp}}
\newcommand{\eex}{\end{exemp}}

\newcommand{\bexs}{\begin{exemps}}
\newcommand{\eexs}{\end{exemps}}

\newcommand{\bprop}{\begin{proposition}}
\newcommand{\eprop}{\end{proposition}}

\newcommand{\bthm}{\begin{theoreme}}
\newcommand{\ethm}{\end{theoreme}}

\newcommand{\bcor}{\begin{corollaire}}
\newcommand{\ecor}{\end{corollaire}}

\newcommand{\blem}{\begin{lemme}}
\newcommand{\elem}{\end{lemme}}

\newcommand{\beqna}{\begin{eqnarray}}
\newcommand{\eeqna}{\end{eqnarray}}

\newcommand{\beqnas}{\begin{eqnarray*}}
\newcommand{\eeqnas}{\end{eqnarray*}}

\newcommand{\cA}{{\mathcal A}}

\newcommand{\LL}{{\rm L}}

\def\I{1\!{\rm l}}
\def\tr{\textmd{trace}\,}

\def\det{{ \rm{det}}}  

\def\Id{{\rm{Id}}} 

\def\diag{{\rm diag}}

\def\cA{{\mathcal A }}

\def\cC{{\mathcal C}}
\def\cD{{\mathcal D}}

\def\cH{{\mathcal  H}}

\def\cL{{\mathcal L }}

\def\cN{{\mathcal N }}

\def\cV{{\mathcal V}}

\def\cA{{\mathcal A }}

\def\cC{{\mathcal C}}
\def\cD{{\mathcal D}}

\def\cH{{\mathcal  H}}

\def\cL{{\mathcal L }}

\def\cN{{\mathcal N }}

\def\cV{{\mathcal V}}

\def\bbC{{\mathbb{C}}}

\newcommand{\bbR}{{\mathbb {R}}}
\newcommand{\bbS}{{\mathbb {S}}} 

\def\un{{\mathbf{1}}}
\def\bfA{{\mathbf{A}}}
\def\bfa{{\mathbf{a}}}
\def\bfB{{\mathbf{B}}}
\def\bfD{{\mathbf{D}}}

\def\bfH{{\mathbf{H}}}

\def\bfM{{\mathbf{M}}}
\def\bfm{{\mathbf{m}}}
\def\bfN{{\mathbf{N}}}
\def\bfO{{\mathbf{O}}}
\def\bfo{{\mathbf{o}}}
\def\bfP{{\mathbf{P}}}
\def\bfS{{\mathbf{S}}}
\def\bfU{{\mathbf{U}}}
\def\bfu{{\mathbf{u}}}
\def\bfV{{\mathbf{V}}}
\def\bfW{{\mathbf{W}}}
\def\bfw{{\mathbf{w}}}
\def\bfX{{\mathbf{X}}}
\def\bfx{{\mathbf{x}}}
\def\bfY{{\mathbf{Y}}}
\def\bfZ{{\mathbf{Z}}}
\def\bfz{{\mathbf{z}}}
\def\bf\Sigma{{\mathbf{\Sigma}}}

\newtheorem{theoreme}{Theorem}[section]
\newtheorem{lemme}[theoreme]{Lemma}
\newtheorem{definition}[theoreme]{Definition}
\newtheorem{proposition}[theoreme]{Proposition}
\newtheorem{corollaire}[theoreme]{Corollary}

\newenvironment{exemp}{\noindent{\bf Example. --- }}{\par}
\newenvironment{exemps}{\noindent{\bf Examples}\benum}{\eenum\par}
\newtheorem{rmq}[theoreme]{Remark}
\newtheorem{rmqs}[theoreme]{Remarks}

\newenvironment{preuve}{\noindent{\it Proof. --- }}
{\hfill\rule{1.3mm}{2mm}\par}
\newenvironment{application}{\noindent{\bf Application. --- }}{\par}
\newenvironment{applications}{\noindent{\bf Applications. ---
}\benum}{\eenum\par}

\theoremstyle{definition}

\title{Matrix Dirichlet processes}

\author{Songzi Li\thanks{songzi.li@bnu.edu.cn, School of Mathematical Sciences, Beijing Normal University, 100875 Beijing, China; Institut de Math\'ematiques de Toulouse, Universit\'e Paul Sabatier,
31062 Toulouse, France;
School of Mathematical Sciences, Fudan University, 200433 Shanghai, China; 
Research supported in part by China Scholarship Council.}}
\date{\today}

\begin{document}
\maketitle

\begin{center}
\begin{minipage}{140mm}
\begin{center}{\textbf{Abstract}}\end{center}  Matrix Dirichlet processes, in reference to their reversible measure, appear in a natural way in many different  models in probability. Applying the language of diffusion operators and the method of boundary equations, we describe Dirichlet processes on the simplex and provide two models of matrix Dirichlet processes, which can be realized by various projections,  through the Brownian motion on the special unitary group, the polar decomposition of complex matrices and also through  Wishart processes. 
\end{minipage}
\end{center}
{\textbf{Key words} :}\ \ Matrix Dirichlet processes, diffusion operators,  polar decomposition, Wishart processes.

\section{Introduction}
The  complex matrix  simplex $\Delta_{n, d}$ is the set of the sequences $(Z^{(1)}, \cdots, Z^{(n)})$ of non negative $d\times d$ Hermitian matrices such that $\sum_{i=1}^{n} Z^{(i)} \leq \Id$, where the inequality is understood in the sense of Hermitian matrices. On the matrix simplex, there exist natural probability measures , with densities $C\det(Z^{(1)})^{a_1-1}\cdots \det(Z^{(n)})^{a_n-1}\det(\Id - \sum^{n}_{i =1}Z^{(i)})^{a_{n+1} -1}$ (see Section \ref{sec.matrix.Dirichlet.complex}). As the natural extensions of the Dirichlet measures on the simplex, they are called matrix Dirichlet measures . 

It turns out that on matrix simplex there exist many diffusion processes which admit matrix Dirichlet measures as reversible ones, and their generators may be diagonalized by a sequence of orthogonal polynomials whose variables are the entries of the matrices. Therefore the matrix simplex appears to be a polynomial domain as described in \cite{BOZ}, see Section \ref{sec.polynomial.models}. 

The purpose of this paper is to describe these diffusion processes that we call matrix Dirichlet processes, referring to their reversible measure.  They appear in a natural way in many different  models in probability : in the polar decomposition of Brownian matrices, in the projection of Brownian motions on $SU(N)$, and also in the projection of Wishart matrices, as we will see below.

 To begin with, we deal with diffusion processes on the simplex which have Dirichlet measures as reversible ones. Dirichlet measures are multivariate generalizations of beta distributions, and play an important role in statistics, such as prior distributions in Bayesian statistics, machine learning, natural language processing, etc. They also appear, together with their associated diffusions processes,  in population biology, for example Wright-Fisher models.   In this paper,  we talk about  Dirichlet processes, by which we mean that these  diffusion processes on the simplex are polynomial models with Dirichlet measures as their reversible measures. We should point out that they have nothing to do with the Dirichlet processes introduced by Ferguson \cite{ferguson}, which have been widely used in statistics, or the Dirichlet processes introduced by F\"ollmer \cite{follmer},  which are related to Dirichlet form.

The matrix Dirichlet measures, as an analogy of Dirichlet distributions, were first introduced by Gupta and Richards~\cite{GuptaRichards}, as  special cases of matrix Liouville measures. They have been deeply studied,  for example by Olkin and Rubin~\cite{olkin}, Gupta and Nagar~\cite{guptanagar}, see also  Letac~\cite{Letac12, LetacMassam} and the references therein. They not only provide  models for multiple random matrices, which are related with orthogonal polynomials, integration formulas, etc., but also reflect the geometry of spaces of matrices, see \cite{guptanagar, Hua}. Therefore, it is natural to consider their corresponding diffusion processes.
It is worth to mention that matrix Jacobi processes, which can be considered as a one matrix case of matrix Dirichlet processes,  were introduced by Doumerc in \cite{Doumerc},  and  studied by Demni in \cite{demni2, demni1}. 

Our interest of this topic not only lies in its importances in statistics and random matrices, but also in the fact that it provides a polynomial model of multiple matrices. In a polynomial model, the diffusion operator (the generator) can be diagonalized by  orthogonal polynomials,  and this leads to the algebraic description of the boundary see Section~\ref{sec.polynomial.models}.  Such polynomial models are quite rare : up to affine transformation, there are 3 models in $\bbR$~\cite{bakrymazet} and 11 models on compact domains in $\bbR^{2}$~\cite{BOZ}. More recently, Bakry and Bressaud~\cite{BakryB} provided new models in dimension 2 and dimension 3, by relaxing the hypothesis in \cite{BOZ} that polynomials are ranked with respect to their natural degree, and investigating the finite groups of $O(3)$ and their invariant polynomials. 

As we will see in this paper, the simplex and the matrix simplex are both polynomial domains, see Section~\ref{sec.polynomial.models}.  Moreover, there exist many different polynomial models on these domians,  which is a quite rare situation : in general, there is just one polynomial diffusion process on a given polynomial domain, up to scaling. The situation here is quite complicated, since we are dealing with a family of matrices. By applying the theory of boundary equations, as introduced in \cite{BakryZ}, we are able to describe  Dirichlet processes on the simplex, and we provide  two models of such matrix Dirichlet processes.  Our two models are found to be realized, via various projections,  through the Brownian motion on the special unitary group, the polar decomposition of complex matrices and also through the Wishart processes. This leads to some efficient ways  to describe image measures in such projections.

This paper is organized as follows. In Section \ref{sec.polynomial.models}, we present the basics on diffusion operators and polynomial models that we will use in this paper. In Section \ref{sec.scalar.simplex}, we introduce the various Dirichlet processes on the simplex, and  describe their  realizations in the special cases  from  spherical Laplacian  or other operators on the sphere, and also from Ornstein-Uhlenbeck or Laguerre  processes. In Section \ref{sec.matrix.Dirichlet.complex}, we introduce our main results, i.e. the description of two polynomial diffusion models, and their realizations through the projections from the Brownian motion on  the special unitary group and  Wishart processes, which are generalizations of the similar results in the scalar (i.e. the simplex) case . Finally in Appendix \ref{polar}, we provide the details of the computations on the polar decomposition of the Brownian motion on complex matrices, which are quite technical.

\section{Symmetric diffusion operators and polynomial models \label{sec.polynomial.models}}
We present in this section a brief introduction to symmetric diffusion processes  and operators, in particular to those diffusion operators which may be diagonalized in a basis of orthogonal polynomials.

Symmetric  diffusion operators are described in~\cite{BGL}, which we refer the reader to for further details. Moreover, for  the particular case of those diffusion associated with orthogonal polynomials, we refer to the paper~\cite{BOZ}. Although the description that we provide below is quite similar to that in ~\cite{BakryB}, we choose to present it here for completeness.

Diffusion processes are Markov processes with continuous trajectories in some open set of $\bbR^n$ or on some manifold, usually given as solutions of stochastic differential equations.  They are described by their infinitesimal generators, which are called diffusion operators.

Diffusion operators are second order differential operators with no zero order terms.
When those operators  have smooth coefficients, they are given in some open subset $\Omega$ of $\bbR^d$ by their action on smooth, compactly supported function $f$ on $\Omega$, 
\beq\label{diffusion.gal} \LL(f) = \sum_{ij} g^{ij}(x) \partial^2_{ij} f+ \sum_i b^i(x)\partial_i f,\eeq
where the symmetric matrix $(g^{ij})(x)$ is everywhere non negative, i.e. the operator $\LL$ is semi-elliptic.

A Markov process $(\xi_t)$ is associated to such a diffusion operator through the requirement that  the process $f(\xi_t)-\int_0^t \LL(f)(\xi_s) ds$ is a local martingale for any function $f$ in the domain of the operator $\LL$.

In this paper we will concentrate on the elliptic case (that is when the matrix $(g^{ij})$ is everywhere non degenerate), and  on the case where this operator is symmetric with respect to some probability measure $\mu$ :  for any smooth functions $f,g$, compactly supported in $\Omega$, we have
\beq\label{IPP} \int_\Omega f \LL(g) d\mu = \int_\Omega g \LL(f) d\mu.\eeq
We say that $\mu$ is a reversible measure for $\LL$ when the associated stochastic process $(\xi_t)$ has a law which is invariant under time reversal, provided that the law at time $0$ of the process is $\mu$. In particular, the measure is invariant : when the associated process $(\xi_t)$ is such that the law of $\xi_0$ is $\mu$, the law of $\xi_t$ is $\mu$ for any time $t>0$.

When $\mu$ has a smooth positive density $\rho$ with respect to the Lebesgue measure, the symmetry property \eqref{IPP} is equivalent to 
\beq\label{eq.density} b^i (x)= \sum_j \partial_j g^{ij}(x) + \sum_j g^{ij} \partial_j \log \rho,\eeq where $b^i(x)$ is the drift coefficient appearing in equation~\eqref{diffusion.gal}. In general,  equation~\eqref{eq.density} allows  to completely determine $\mu$  from the data in $\LL$, up to some normalizing constant .

Now we introduce the carré du champ operator $\Gamma$. Suppose that we have some dense algebra $\cA$ of functions  in $\cL^2(\mu)$ which is stable under the operator $\LL$ and contains the constant functions. Then for $(f,g)\in \cA$  we define
\beq\label{def.Gamma} \Gamma(f,g) = \frac{1}{2}( \LL(fg)-f\LL(g)-g\LL(f)).\eeq

If $\LL$ is given by equation~\eqref{diffusion.gal}, and the elements of $\cA$ are at least $\cC^2$,   we have
$$\Gamma(f,g)=\sum_{ij} g^{ij} \partial_if \partial_j g,$$ so that $\Gamma$ describes in fact the second order part of $\LL$.  The semi-ellipticity of $\LL$ gives rise to the fact that $\Gamma(f,f)\geq 0$, for any $f\in \cA$.

Applying formula~\eqref{IPP} with $g=1$, we obtain $\int_\Omega \LL f d\mu=0$ for any $f\in \cA$. Then with ~\eqref{IPP} again, we see immediately that for any $(f,g)\in \cA$
\beq\label{IPP2} \int_\Omega f\LL(g) d\mu = -\int_\Omega \Gamma(f,g) d\mu,\eeq so that the knowledge of $\Gamma$ and $\mu$ describes entirely the operator $\LL$. Such a triple  $(\Omega, \Gamma, \mu)$ is called a Markov triple in~\cite{BGL}.

By~\eqref{diffusion.gal}, we see that $\LL(x^i)= b^i$ and $\Gamma(x^i, x^j)= g^{ij}$. The operator $\Gamma$ is called the co-metric, and in our system of coordinates is described by a matrix $\Gamma= \big(\Gamma(x^i, x^j)\big)= (g^{ij})$.

In our setting, we will  always  assume that $\Omega$ is bounded and choose $\cA$ to be the set of polynomials. Since polynomials are not compactly supported in $\Omega$, the validity of equation~\eqref{IPP} requires extra conditions on the coefficients $(g^{ij})$ at the boundary of $\Omega$, which we will describe below. 

The fact that $\LL$ is a second order differential operator implies the  change of variable formulas. Whenever $f=(f_1, \cdots, f_n)\in \cA^n$ and $\Phi(f_1, \cdots, f_n)\in \cA$, for some smooth function $\Phi  : \bbR^n\mapsto \bbR$, we have
\beq\label{chain.rule.L}\LL(\Phi(f))= \sum_i \partial_i \Phi(f) \LL(f_i) + \sum_{ij} \partial^2_{ij} \Phi(f) \Gamma(f_i,f_j)\eeq
and also
\beq\label{chain.rule.Gamma} \Gamma(\Phi(f),g) = \sum_i \partial_i \Phi(f) \Gamma(f_i,g).\eeq
When $\cA$ is the algebra of polynomials, properties~\eqref{chain.rule.L} and~\eqref{chain.rule.Gamma} are equivalent.


 An important feature in the examples described in this paper is the notion of image. Whenever we have a diffusion operator $\LL$ on some set $\Omega$, it may happen that we find some functions $(X_1, \cdots, X_p)$ in the algebra $\cA$ such that $\LL(X_i)= B^i(X)$ and $\Gamma(X_i, X_j)= G^{ij}(X)$ where $X= (X_1, \cdots, X_p)$. Then we say that we have a closed system.  If $(\xi_i)$ is the Markov diffusion process with generator $\LL$, $(\zeta_t)= X(\xi_t)$ is again a diffusion process, with its generator expressed in coordinates $(X_1, \cdots, X_p)$
 $$\hat \LL = \sum_{ij} G^{ij}(X) \partial^2_{ij} + B^i(X)\partial_i.$$
 Moreover, when $\LL$ has a reversible  probability measure $\mu$, $\hat \LL$ has the image measure of $\mu$ through the map $X$ as its reversible measure  . Thanks to equation~\eqref{eq.density}, this is often an efficient way to determine image measure, which will be used many times in this paper.

As mentioned above, we will restrict our attention to the  elliptic case. Here we expect $\LL$ to have a self adjoint extension (not unique in general), thus it has a spectral decomposition. Also 
we expect that the spectrum is discrete,  thus that  it has eigenvectors, which we will require to be polynomials in the variables $(x^i)$.  Moreover, we will require that those polynomial eigenvectors to be ranked according to their degrees, i.e., if we denote by  $\cH_n$  the space of polynomials with total degree at most $n$, then for each $n$ we need that there exists an orthonormal basis of $\cH_n$ which is made of eigenvectors for $\LL$. Equivalently, we require that $\LL$ maps $\cH_n$ into itself. This situation is quite rare, and imposes some strong restriction on the domain $\Omega$ that we will describe below.

When we have such a Markov triple  $(\Omega, \Gamma, \mu)$, where  $\Omega$ has a piecewise smooth (at least $\cC^1$) boundary, $\mu$ has a smooth density with respect to the Lebesgue measure on $\Omega$. Then we call $\Omega$ a polynomial domain and $(\Omega, \Gamma, \mu)$ a polynomial model.

In dimension 1 for example, up to affine transformations, there are only 3 cases of polynomial models, corresponding to the Jacobi, Laguerre and Hermite polynomials, see for example~\cite{bakrymazet}.

\begin{enumerate}
\item \label{Hermite.case} The Hermite case corresponds to the case where  $\Omega= \bbR$, $\mu$ is the Gaussian measure  $\frac{e^{-x^2/2}}{\sqrt{2\pi}}\,dx$ on $\bbR$ and $\LL$ is  the  Ornstein--Uhlenbeck operator
\beq\label{def.OU}\LL_{OU}=\frac{d^2}{dx^2}-x\frac{d}{dx}.\eeq  The Hermite polynomial $H_n$ of degree $n$ satisfy $\LL_{OU} P_n= -nP_n$.

\item\label{Laguerre.case} The Laguerre polynomials correspond to the case where $\Omega= (0, \infty)$, the measure $\mu$  depends on a parameter $a>0$ and is  $\mu_a(dx)=C_ax^{a-1}e^{-x}\,dx$ on $(0, \infty)$,
and $\LL$ is  the Laguerre operator

\beq\label{def.Laguerre}\LL_a=x\frac{d^2}{dx^2}+ (a-x)\frac{d}{dx}.\eeq The Laguerre  polynomial $L^{(a)}_n$ with  degree $n$ satisfies $\LL_a L^{(a)}_n= -nL_n^{(a)}$.

\item \label{Jacobi.case}The Jacobi polynomials  correspond to the case where $\Omega= (-1,1)$,  the measure $\mu$ depends on two parameters $a$ and $b$, $a,b>0$ and is  is  $\mu_{a,b}(dx)=C_{a,b}(1-x)^{a-1}(1+x)^{b-1}\,dx$ on $(-1,1)$, and $\LL$ is  the Jacobi operator
\beq\label{def.Jacobi}\LL_{a,b}=(1-x^2)\frac{d^2}{dx^2}- \big(a-b+(a+b)x\big)\frac{d}{dx}.\eeq
The Jacobi polynomial  $(J^{(a,b)}_n)_n$ with degree $n$  satisfy $$\LL_{a,b}J_n^{(a,b)}= - n(n+a+b-1)J_n^{(a,b)}.$$
\end{enumerate}

For the general case, when $\Omega$ is bounded, we recall here some results in~\cite{BOZ}.

\bprop\label{prop.poly.model.gal} Let $(\Omega, \Gamma, \mu)$  be a polynomial model in $\bbR^d$. Then,  with $\LL$ described by equation~ \eqref{diffusion.gal}, we have
\benum

\item\label{prop1.pnt1} For $i= 1, \cdots, d$, $b^i$ is a polynomial with $\deg(b^i)\leq 1$.

\item \label{prop1.pnt2} For $ i, j= 1, \cdots, d$, $g^{ij}$ is a polynomial with $\deg(g^{ij})\leq 2$.

\item \label{prop1.pnt3} The boundary $\partial\Omega$ is included in the algebraic set $\{\det(g^{ij})=0\}$.

\item \label{prop1.pnt4}  If $\{P_1\cdots P_k=0\}$ is the reduced equation of the boundary $\partial\Omega$ (see remark~\ref{rmk.reduced.bdry.eq} below), then, for each $q= 1, \cdots k$, each $i= 1, \cdots d$, one has
\beq\label{boundary.eq}\Gamma( \log P_q, x_i)= L_{i,q},\eeq where
$L_{i,q}$ is a polynomial with $\deg(L_{i,q})\leq 1$;

\item\label{prop1.pnt5}  All the measures $\mu_{\alpha_1, \cdots, \alpha_k}$ with densities $C_{\alpha_1, \cdots, \alpha_k}|P_1|^{\alpha_1}\cdots |P_k|^{\alpha_k}$ on $\Omega$, where the $\alpha_i$ are such that the  density is  is integrable on $\Omega$, are such that $(\Omega, \Gamma, \mu_{\alpha_1, \cdots, \alpha_k})$ is a polynomial model.

\item\label{prop1.pnt6}  When the degree of $P_1\cdots P_k$ is equal to the degree of $\det(g^{ij})$, there are no other measures.

\eenum

Conversely, assume that some bounded domain $\Omega$ is such that the boundary $\partial\Omega$ is included in an algebraic surface and has reduced equation $\{P_1\cdots P_k=0\}$. Assume moreover that there exists a matrix $(g^{ij}(x))$ which is positive definite in $\Omega$ and such that each component $g^{ij}(x)$ is a polynomial with degree at most $2$.  Let $\Gamma$ denote the associated  carr\'e du champ operator. Assume moreover that equation~\eqref{boundary.eq} is satisfied for any $i= 1, \cdots, d$ and any $q= 1,\cdots, k$, with $L_{i,q}$ a polynomial with degree at most $1$.

 Let
$(\alpha_1, \cdots, \alpha_k)$ be such  that  the $|P_1|^{\alpha_1}\cdots |P_k|^{\alpha_k}$ is integrable  on $\Omega$ with respect to the Lebesgue measure, and denote 
$$\mu_{\alpha_1, \cdots, \alpha_k} (dx)= C_{\alpha_1, \cdots, \alpha_k}|P_1|^{\alpha_1}\cdots |P_k|^{\alpha_k}dx,$$ where $C_{\alpha_1, \cdots, \alpha_k}$ is the normalizing constant such that $\mu_{\alpha_1, \cdots, \alpha_k}$ is a probability measure.

 Then  $(\Omega, \Gamma,\mu_{\alpha_1, \cdots, \alpha_k})$ is a polynomial model.

\eprop

\brmq\label{rmk.reduced.bdry.eq} We say that $\{P_1\cdots P_k=0\}$ is the reduced equation of the boundary $\partial\Omega$ when
\benum
\item The polynomials $P_i$ are not proportional to each other.
\item For $i= 1, \cdots k$, $P_i$ is an irreducible polynomial, both in the real and the complex field.

\item For each $i= 1, \cdots, k$, there exists at least one regular point of the boundary $\partial\Omega$ such that $P_i(x)= 0$.

\item For each regular point $x\in \partial\Omega$, there exist a neighborhood $\cV$  and of $x$  and some $i$ such that
$\partial\Omega\cap \cV= \{P_i(x)=0\}\cap \cV$.

\eenum

In particular, this does not mean that any point satisfying $P_i(x)=0$ for some $i$ belongs to $\partial\Omega$.

\ermq

\brmq The determination of the polynomial domains therefore amounts to the determination of the domains $\Omega$ with an algebraic boundary, with the property that  the reduced equation of $\partial\Omega$ is such that the set of equations~\eqref{boundary.eq} has a non trivial solution, for $g^{ij}$ and $L_{i,q}$. Given the reduced equation of $\partial\Omega$,  equations \eqref{boundary.eq} are linear homogeneous ones in the coefficients of the polynomials $g^{ij}$ and of the polynomials $L_{i,k}$. Unfortunately, in general we need much more equations to determine the unknowns uniquely, and this requires very strong constraints on the polynomials appearing in the reduced equation of the boundary.  We will see that both the simplex and the matrix simplex (in the complex and real case) are such domains where the choice of the co-metric $\Gamma$ is not unique.

\ermq

\brmq The set of equations~\eqref{boundary.eq}, which are central in the study of polynomial models, may be reduced to less equations, when $k>1$.  Indeed, if we set $P= P_1\cdots P_k$, it reduces to
\beq\label{boundary.eq2}\Gamma(x_i,\log P)= L_i, \quad \deg(L_i)\leq 1.\eeq 
In fact assume that this last  equation holds with some polynomial $L_i$,   then on the regular part of the boundary described by $\{P_q(x)=0\}$,  we have $\Gamma(x_i,  P_q)= 0$ since
$$\Gamma(x_i, P_q)= P_q(L_i-\sum_{l\neq q} \frac{\Gamma(x_i,P_l)}{P_l}).$$
Therefore, $\Gamma(x_i,P_q)$ vanishes on the regular part, and $P_q$ being irreducible, it divides $\Gamma(x_i, P_q)$. This leads to $\Gamma(x_i, P_q)= L_{i,q} P_q$, where $\deg(L_{i,q})\leq 1$. Thus we obtain the equation~\eqref{boundary.eq}.

\ermq

\brmq\label{rmq.def.polynomial.domain} A bounded  polynomial domain is therefore any  bounded domain $\Omega$ with algebraic boundary, on which there exists a symmetric matrix $(g^{ij})$  with entries which are polynomials with degree at most 2, that is positive definite on $\Omega$ and defines on $\Omega$ an operator $\Gamma$  satisfying equation~\eqref{boundary.eq2}. \ermq

 In~\cite{BOZ}, a complete description of all polynomial domains and models in dimension 2 is provided : 11 different cases are given up to affine transformations. This description only relies on algebraic considerations on those algebraic curves in the plane where the boundary condition~\eqref{boundary.eq} has a non trivial solution. This reflects the fundamental role played by the boundary equation.

 Among the 11 bounded domains provided by the classification in~\cite{BOZ}, the triangle appears to be one of the few ones (with the unit ball and a particular case of the parabolic bi-angle) where the metric is not unique. In higher dimension, the triangle  may be generalized to  the simplex, which appear to be also a polynomial domain.  Here we extend it furthermore  to the matrix simplex.

\section{Dirichlet measure on the simplex\label{sec.scalar.simplex}}

In this section, we recall a few facts about the simplex and provide some diffusion operators on it which shows that the simplex is a polynomial domain in the sense of Section~\ref{sec.polynomial.models}.

\bdefi The $n$ dimensional simplex $\Delta_{n}\subset \bbR^{n}$ is the set

$$
\Delta_{n} = \{0\leq x_{i} \leq 1,  ~ i= 1,\cdots n, \sum^{n}_{i = 1}x_{i} \leq 1,  \}.
$$
Given  $\bfa = (a_1, a_2,..., a_{n+1}) \in \mathbb{R}^{n+1}$, where $a_{i}>0$, $ i =  1,..., n+1$,  the Dirichlet distribution $D_{\bfa}$
is the probability measure on $\Delta_n$  given by
$$
\frac{1}{B_{\bfa}}x_{1}^{a_1 - 1}...x_{n}^{a_n - 1}(1 - x_1- ...- x_{n})^{a_{n+1} - 1}\mathbf{1}_{\Delta_{d}}(x_1, ..., x_n)dx_1...dx_{n},
$$
where $B_{\bfa} = \frac{\Gamma(a_1)...\Gamma(a_{n+1})}{\Gamma(a_1 + ...+a_{n +1})}$ is the normalizing constant.
\edefi
The Dirichlet measure can be considered as a $n$-dimensional  generalization of the beta distribution on the real line,
$$
\beta(a_1, a_2) = \frac{1}{B(a_1, a_2)}x^{a_1 -1}(1 - x)^{a_2 - 1}\mathbf{1}_{(0, 1)}(x)dx,
$$
which is indeed $D_{(a_1, a_2)}$.

It turns out that the simplex $\Delta_n$ is a polynomial diffusion domain in $\bbR^n$ in the sense of \cite{BOZ}, as described in Section~\ref{sec.polynomial.models}.  More precisely, there exist many different polynomial models on the simplex which admit the Dirichlet measure as their reversible measure.

\bthm\label{th.Dirichlet.simplex}
Let $\mathbf{A}=(A_{ij})_{i,j=1\cdots ,n+1}$ be a  symmetric $(n+1)\times (n+1)$ matrix    where all the coefficients $A_{ij}$ are non negative, and $A_{ii}=0$. Let $\mathbf{a}= (a_1, \cdots, a_{n+1})$ be a $(n+1)$-tuple of positive real numbers.  Let $\LL_{\bfA, \bfa}$ be the symmetric diffusion operator  defined on the   simplex $\Delta_n$ defined by
\beqna\label{gamma.dirichlet.simplex}
\Gamma_{\bfA}(x_i, x_j) &=& -A_{ij}x_ix_j + \delta_{ij}\sum^{n+1}_{k=1}A_{ik}x_kx_i ,\\
 \LL_{\bfA,\bfa}(x_i) &=& -x_i\sum^{n+1}_{j=1} A_{ij}a_j + a_i\sum^{d+1}_{j=1}A_{ij}x_j,
\eeqna
where $x_{n+1} = 1 - \sum^{n}_{j=1}x_j$.

Then, as soon as $A_{ij}>0$ for all $i,j= 1,\cdots, n+1$, $i\neq j$, $(\Delta_n, D_{\bfa}, \Gamma_{\bfA})$ is a polynomial model, and the operator $\LL_{\bfA,\bfa}$ is elliptic on $\Delta_n$. Moreover, any  polynomial diffusion model on $\Delta_n$ having $D_{\bfa}$ as its reversible measure  has this form (however without the requirement that the coefficients $A_{ij}$ are positive). 
\ethm

Observe that the condition $A_{ii}=0$ is irrelevant in formulas \eqref{gamma.dirichlet.simplex}, since the coefficient $A_{ii}$ vanishes in the formulas.

\bpf

Let us prove first that any polynomial model with the usual degree has this form.
 According to Proposition~\ref{prop.poly.model.gal}, to be a polynomial model,  $\Gamma(x_i, x_j)$ must be a polynomial no more than degree 2 and satisfy the boundary equation \eqref{boundary.eq} :  for $1 \leq i, j \leq n$,
\beqna
\label{eq1}\Gamma_{\bfA}(x_i, \log x_{j}) &=& L^{i,j},\\
\label{eq2}\Gamma_{\bfA}(x_i, \log x_{n+1}) &=& L^{i,n+1},
\eeqna
where $\{L^{i,j}, 1 \leq j \leq n+1\}$  are polynomials with degree at most  1.

From equation~\eqref{eq1}, we get
$$
\Gamma_{\bfA}(x_i, x_j) = x_jx_jL^{ij},
$$
so that $\Gamma_{\bfA}(x_i,x_j)$ is divisible by $x_i$, and similarly by $x_j$. Therefore when $i\neq j$, there exists a constant $A_{ij}$, with $A_{ij}= A_{ij}$, such that
$$
\Gamma_{\bfA}(x_i, x_j) = -A_{ij}x_jx_i.$$

When $i=j$, $\Gamma_{\bfA}(x_i,x_i)$ is divisible by $x_i$ and
from equation~\eqref{eq2},  we obtain
$$
\sum^{n}_{j=1}\Gamma_{\bfA}(x_i, x_j) = -x_{n+1}L^{i,n+1},
$$
which writes
$$
\Gamma_{\bfA}(x_i,x_i)= x_i(\sum_{j=1, j\neq i}^{n}A_{ij}x_j) -x_{n+1}L^{i,n+1}.
$$
This implies that $x_i$ divides $L^{i,n+1}$, and therefore that $L^{i,n+1}= -A_{i(n+1)} x_i$, so that
$$\Gamma_{\bfA}(x_i,x_i)= x_i \sum_{i=1, j\neq i}^{n+1} A_{ij} x_j.$$

Conversely, it is quite immediate every operator $\Gamma_\bfA$ defined by equation~\eqref{eq1} satisfies the boundary condition on the simplex.

The formulas for $\LL_{\bfA, \bfa}$ are just a direct consequence of the reversible measure equation~\eqref{eq.density}.

The ellipticity of $\Gamma_{\bfA}$ when all the coefficients $A_{ij}, i\neq j$ are positive  is a particular case of the real version  of Theorem~\ref{matrix.dirichlet}, which will be proved in Section~\ref{sec.matrix.Dirichlet.complex}.

\epf

When all the coefficients $A_{ij}$ are equal to 1, and when all the parameters $a_i$ are half integers, the operator \eqref{gamma.dirichlet.simplex} is an image of the spherical Laplacian.  Indeed,   consider the unit sphere $\sum_{i= 1}^{N} y_i^2=1$  in $\bbR^N$, and split the set $\{1, \cdots, N\}$ in a partition of $n+1$ disjoint subsets $I_1, \cdots, I_{n+1}$ with respective size $p_1, \cdots, p_{n+1}$. For $i= 1, \cdots, n+1$, set $x_i= \sum_{j\in I_i} y_j^2$, so that $x_{n+1}= 1-\sum_{i=1}^{n} x_i$.   Then it is easy to check that, for the spherical Laplace operator $\Delta_{\bbS^{N-1}}$ in $\bbR^N$ (see \cite{BGL} for details), $\Delta_{\bbS^{N-1}} (x_i)$ and $\Gamma_{\bbS^{N-1}}(x_i,x_j)$ coincide with those given in equations \eqref{gamma.dirichlet.simplex} whenever $a_i = \frac{p_i}{2}$.

The next proposition generalizes this geometric interpretation for the general  choice of the parameters $A_{ij}$.

\bprop\label{prop.general.scalar.case}
Let $I_1, \cdots, I_{n+1}$ be a partition of $\{1, \cdots, N\}$ into disjoint sets with size $p_1, \cdots, p_{n+1}$.  For $i, j\in \{1, \cdots, n+1\}$, with $i\neq j$,  let $\LL_{i,j}  $ be the operator acting on the unit sphere $\bbS^{N-1}\subset \bbR^N$ as
$$\LL_{i,j} = \sum_{p\in I_i, q\in I_j} (y_p\partial_{y_q}-y_q\partial_{y_p})^2.$$

Then, setting $x_i= \sum_{p\in I_i} y_i^2$, $\{x_1, \cdots, x_{n+1}\}$ form a closed system for any $\LL_{i,j}$, and the image of
$$\frac{1}{2}\sum_{i < j} A_{ij} \LL_{i,j}$$ is the operator $\LL_{\bfA,\bfa}$, where $a_i= \frac{p_i}{2}$.
\eprop


The proof follows from a direct application of the change of variable formulas~\eqref{chain.rule.L} and~\eqref{chain.rule.Gamma} applied to the functions $\{x_p\}$. For any operator $\LL_{i,j}$,   the associated carré du champ operator is $\Gamma_{i,j}(f,g)= \sum_{p\in I_i, q\in I_j} (y_p\partial_{y_q}
-y_q\partial_{y_p})(f)(y_p\partial_{y_q} -y_q\partial_{y_p})(g)$.  We just  have to identify $\LL_{i,j} (x_p)$ and $\Gamma_{i,j}(x_p,x_q)$ through a direct and easy computation,  and observe that it fits with the coefficients of $A_{ij}$ in the same expression in the definition of $\LL_{\bfA,\bfa}$.

It is worth to observe that the spherical Laplace operator is nothing else than $\sum_{i,j} \LL_{i,j}$. The disappearance of the operators $\LL_{i,i}$ in the general form for $\LL_{\bfA,\bfa}$ comes from  from the fact that the action of $\LL_{i,i}$ on any of the variables $x_p$ vanishes.

To come back to the case where all the coefficients $A_{ij}$ are equal to 1, and as a consequence of the previous observation, the Dirichlet measure on the simplex is an image of the uniform measure on the sphere when the parameters are half integers, through the map that we just described $(y_1, \cdots, y_N)\mapsto (x_1, \cdots, x_{d})$ where $x_i= \sum_{j\in I_i} y_j^2$.   In the same way that the uniform measure on a sphere is an image of a Gaussian measure through the map $y\in \bbR^N\mapsto \frac{y}{\|y\|}\in \bbS^{N-1}$, the Dirichlet measure may be constructed from the random variables $\sum_{j\in I_i} y_j^2$ where $y_j$ are independent standard Gaussian variables. More precisely, if  $(y_1, \cdots, y_N)$ are independent real valued Gaussian variables, setting $S_i= \sum_{j\in I_i} y_j^2$ and $S= \sum_{i=1}^{d+1} S_i$, we see that
$(\frac{S_1}{S},\cdots,  \frac{S_d}{S} )$ follows the Dirichlet law $D_{\bfa}$ whenever $a_i= \frac{p_i}{2}$.

Since the variables $S_i$ have $\gamma$  distribution $\frac{1}{\Gamma(\frac{p_i}{2})}x^{\frac{p_i}{2} - 1}e^{-x/2}dx$, it is not surprising that this representation of Dirichlet laws may be extended for the general case, that is when the parameters $a_i$ are no longer half integers, by replacing norms of Gaussian vectors by independent variables having $\gamma$ distribution.  We quote the following proposition from~\cite{Letac12}, which gives a construction of Dirichlet random variable through gamma distributions $\gamma_{\alpha}$ on $\bbR^{+}$ given by
$$\gamma_{\alpha}(dx) = \frac{1}{\Gamma(\alpha)}x^{\alpha-1}e^{-\frac{x}{\beta}}\mathbf{1}_{(0, \infty)}(x)\, dx,
$$
where $\alpha> 0$ and  $\Gamma(\alpha) = \int^{\infty}_0 x^{\alpha-1}e^{-x}dx$.
\bprop\label{prop6}
Consider independent random variables $x_1,..., x_n, x_{n+1}$ such that each $x_k$ has gamma distribution $\gamma_{\alpha_k, 2}$. Define $S = \sum^{n+1}_{k=1}x_k$. Then $S$ is independent of $\frac{1}{S}(x_1,..., x_n)$, and the distribution of $\frac{1}{S}( x_1,..., x_n)$ has the density $D_{(\alpha_1, ..., \alpha_n)}$.
\eprop

Indeed, Proposition~\ref{prop6} is still valid at the level of the processes. Start first from  an $N$-dimensional Ornstein-Uhlenbeck process, which admits the  standard Gaussian measure in $\bbR^N$ as reversible measure.  One should be careful here, since the spherical Brownian motion is not directly the image of an Ornstein-Uhlenbeck process through  $x\mapsto \frac{x}{\|x\|}$. Indeed,  writing the Ornstein-Uhlenbeck operator $\LL_{OU}= \Delta-x\cdot\nabla$ in  $\bbR^N\setminus\{0\}$ in polar coordinates $(r= \|x\|, \phi= \frac{x}{\|x\|})\in (0, \infty)\times \bbS^{N-1}$,  we derive
$$\LL_{OU}= \partial^2_r + (\frac{N-1}{r} -r) \partial_r + \frac{1}{r^2} \Delta_{\bbS^{N-1}},$$
where $\Delta_{\bbS^{N-1}}$ is the spherical Brownian Laplace operator acting on the variable  $\phi\in \bbS^{N-1}$.

The structure of $\LL_{OU}$ shows that it appears as a skew product (or warped product) of a one dimensional operator with the spherical Laplace operator. In the probabilist interpretation, if $\xi_t$ is an Ornstein-Uhlenbeck process  in $\bbR^N$, and if $\rho_t = \| \xi_t\|$  and $\theta_t= \frac{\xi_t}{\|\xi_t\|}$, then $\rho_t$ is a diffusion process (with generator $\partial^2_r + (\frac{N-1}{r} -r) \partial_r$), and the pair $(\rho_t, \theta_t)$ is also a diffusion process, while it is not the case for $\theta_t$ alone. We have to change the clock at which $\theta_t$ is running,  namely set $\hat \theta_t=  \theta_{C_t}$ where $C_t= \int_0^t \rho_s^2 ds$, such that $\hat \theta_t$ is a spherical Brownian motion which is independent of  $\rho_t$.  We may also consider the process $\theta_t$ alone conditioned on $\rho_t=1$. However, since the reversible law of $\hat \theta_t$ is the uniform measure on the sphere, we also see that, when the  law of $\xi_0$ is  the reversible law of the the process $\xi_t$ (that is the standard $\cN(0,1)$ Gaussian law in $\bbR^N$),  the law of $\theta_t$ given $\rho_t$ does not depend on $\rho_t$ for fixed time $t$, such that the law of $\theta_t$ is the uniform measure on the sphere and $\theta_t$ and $\rho_t$ are independent.

Starting again with some Ornstein-Uhlenbeck process in $\xi_t= (\xi_{t}^{(1)}, \cdots, \xi_t^{(N)})$ in $\bbR^N$, and cutting the set of indices in $n+1$ parts $I_i$ as above, with $|I_i|= p_i$, we may now consider the variables $\sigma_t^{(i)}= S_i(\xi_t)$ and denote $\sigma_t= S(\xi_t)$, where $S_i(x)= \frac{1}{2}\sum_{j\in I_i} x_j^2$ and $S(x)= \sum_i S_i(x)$.  Then $\sigma_t^{(i)}$ are independent Laguerre processes with generator
$\LL_{a_i}(f)= 2(xf''+ (a_i-x)f')$,  where $a_i= \frac{p_i}{2}$.
One sees  that  $\sigma_t$ is again a Laguerre process with generator $\LL_a(f)= 2(xf''+ (a-x)f')$, with $a= N/2$. Then, setting $\zeta^{(i)}_t=  \frac{\sigma^{(i)}_t}{\sigma_t}$ and $\zeta_t= (\zeta^{(1)}_t, \cdots, \zeta^{(n)}_t)$,
$\zeta_t$ takes values in the simplex $\Delta_n$, and a simple computation shows that
$(\sigma_t, \zeta_t)$ is a Markov process with generator $\LL_a + \frac{1}{S}\LL_{\bfA,\bfa}$, where $\LL_{\bfA, \bfa}$ is the operator defined in Theorem \ref{th.Dirichlet.simplex}, with  $A_{ij}= 2$ and $a_i= \frac{p_i}{2}$.

This may be extended to the general case where the parameters $a_i$ are no longer half integers. Starting from the standard Laguerre operator, we have

\bprop \label{simplex.construction1} Let  $\bar a= (a_1, \cdots, a_{n+1})$ be positive integers and $\Sigma_t=(\sigma_t^{(1)}, \cdots, \sigma^{(n+1)}_t) $ be independent Laguerre  processes on $\bbR_+$ with generator
$\LL_{a_i}(f)= xf''+ (a_i-x)f'$. Let $\sigma_t= \sum_{i=1}^{n+1} \sigma_t^{(i)}$ and $\zeta_t= (\frac{\sigma_t^{(1)}}{\sigma_t}, \cdots , \frac{\sigma_t^{(d)}}{\sigma_t})$.  Then, the pair $(\sigma_t, \zeta_t)\in \bbR_+ \times \Delta_{n}$ is a diffusion process with generator
$$S\partial^2_S+ (a-S)\partial_S+ \frac{1}{S}\LL_{\bfA,\bfa}$$ where $\LL_{\bfA,\bfa}$ is defined in Theorem~\ref{th.Dirichlet.simplex} with $A_{ij}=1$ for $i\neq j$, $\bfa = \bar{a}$ and $a= \sum_{i=1}^{n+1} a_i$.
\eprop

Here, if the coordinates that we choose for a variable in $\bbR_+ \times \Delta_n$ are $(S, x)$, where $x= (x_1, \cdots, x_n)\in \Delta_n$,  then the action of the generator is $\LL_{a}$ on the variable $S$ and $\frac{1}{S}\LL_{\bfA,\bfa}$ on the variable $x$. We are in the situation of a warped product, and $\sigma_t$  is a Laguerre diffusion process with generator $\LL_{a}$. Setting $C_t= \int_0^t \sigma_s ds$,  $\zeta_{C_t}$ is a diffusion process on $\Delta_n$ which is independent of $\sigma_t$ and has generator $\LL_{\bfA,\bfa}$.  Moreover, the reversible measure of $\Sigma_t$ is the  product on independent $\gamma_{a_i}$ measures on $\bbR_+$, and its image on $\bbR_+\times \Delta_n$  through the map
$y= (y^{(1)}, \cdots, y^{(n+1)})\in \bbR_+^{n+1} \mapsto (S= \sum_{i=1}^{n+1} y_i, x_i= \frac{y_i}{S})$ is the product of   a $\gamma_{a}$ measure on $\bbR_+$ and of a $D_{\bfa}$ measure on $\Delta_n$. We thus get another proof of Proposition \ref{prop6}.

\bpf

Following the notations of Section~\ref{sec.polynomial.models}, we consider he generator of the process $\Sigma_t$, which is
$$\LL_{\bar a}=\sum_{i= 1}^{n+1} y_i\partial^2_{y_i} f+  \sum_{i= 1}^{n+1} (a_i-y_i)\partial_{y_i} f,$$ where $y= (y_1, \cdots, y_{n+1})\in \bbR_+^{(n+1)}$.  Let $\Gamma_{\bar a}$ be its associated carré du champ operator. With $S= \sum_{i=1}^{n+1} y_i$ and $z_i=\frac{y_i}{S}$, $i= 1, \cdots, n$,
we just have to check that
\benum
\item $\LL_{\bar a}(S)= a-S$;
\item $\LL_{\bar a}(z_i)=\frac{1}{S}\big(a_i - az_i\big) $;
\item $\Gamma_{\bar a}(S,S)= S$;
\item $\Gamma_{\bar a}(S,z_i)=0$;
\item $\Gamma_{\bar a}(z_i,z_j)=\frac{1}{S}(\delta_{ij} z_i-z_iz_j) $;
\eenum
This follows directly from a straightforward computation.

\epf

Proposition~\ref{simplex.construction1} provides a construction of the process on the simplex with generator $\LL_{\bfA, \bfa}$ for the general $\bfa$ but only when $A_{ij}= 1$ for $i\neq j$. To obtain a construction in the general case, we may use the same generalization that we did on the sphere with the operators $\LL_{i,j}= \sum_{k\in I_i, l\in I_j} (y_k\partial_l-y_l\partial_k)^2$, where $y= (y_i)\in \bbS^n$, i.e., consider the operator
$$\LL^{\bfA}_{OU}= \partial^2_r + (\frac{N-1}{r} -r) \partial_r + \frac{1}{2r^2}\sum_{i < j}A_{ij}\LL_{i, j},$$
where we use $\frac{1}{2}\sum_{i < j}A_{ij}\LL_{i, j}$ instead of $\Delta_{\bbS^{N-1}}$ in the Ornstein-Uhlenbeck operator.
Its reversible measure is still $e^{-\frac{1}{2}\|x\|^{2}}dx_{1}...dx_{N+1}$ on $\bbR^{N+1}$. Now let $\LL^{\bfA}_{OU}$ act on the variables $x^{i}= \sum_{p\in I_i  } y_p^2$, and $S = \sum_{i}x_i$ $z_i=\frac{x_i}{S}$, $i= 1, \cdots, n$, then it is easy to check that the image of $\LL^{\bfA}_{OU}$ is $\LL_{\bfA, \bfa}$ with $a_i = \frac{p_i}{2}$. 

In fact this construction is essentially the same as Proposition \ref{prop.general.scalar.case}. Consider a diffusion process $\bfx = (x_{1}, \cdots ,x_{n+1})$ on $\bbR^{n+1}$ with a warped product symmetric diffusion generator
$$
\LL = \LL_{r} + f(r)\Delta_{\bbS^{n}},
$$
where $\LL_{r}$ is a symmetric diffusion operator on $r_t = \|x_{t}\|$, $f: \bbR^{+} \rightarrow \bbR^{+}$ is a function of $r$. Define $y_{i} = x^{2}_{i}$ and $z_{i} = \frac{y_{i}}{r^{2}}$ for $1 \leq i \leq (n+1)$, then the generators of $\bfz = (z_{1}, \cdots ,z_{n+1})$ are indeed the image of $f(r)\Delta_{\bbS^{n}}$ through the map from $\bfx$ to $\bfz$. Following the same procedure as in Theorem \ref{th.Dirichlet.simplex}, we can always construct the Dirichlet process with a general parameter~$\bfA$.

\section{Matrix Dirichlet Processes \label{sec.matrix.Dirichlet.complex}}

\subsection{The complex matrix simplex and their associated Dirichlet measures\label{subsec.matrix.Dirichlet.domain}}
We now introduce the matrix simplex, the matrix Dirichlet processes and  its corresponding diffusion operators.  In particular, we give two families of models. The first one comes from  extracted matrices from the Brownian motion on $SU(N)$,  and is an analogy to the construction of the scalar Dirichlet process from spherical Brownian motion.  Similarly, Proposition~\ref{simplex.construction1}  may be extended to a similar construction from Wishart matrix processes, which are a natural matrix extension of Laguerre processes. Another construction arises naturally when one studies the polar decomposition of Brownian matrices, in which case the spectral projector provides a degenerate form of such matrix Dirichlet processes.

Considering the real version of matrix Brownian motions, it is not hard to derive the real counterparts of our results, that we do not develop here.

First we define the complex matrix generalization  $\Delta_{n, d}$ of the $n$-dimensional simplex  as the set of $n$-tuple of Hermitian  non negative matrices $\{\bfZ^{(1)}, \cdots, \bfZ^{(n)}\}$ such that
\beq\label{matrix.simplex}
\sum_{k=1}^n \bfZ^{(k)} \leq \Id.
\eeq

Accordingly, the complex matrix Dirichlet measure $D_{\bfa}$ on $\Delta_{n,d}$ is given by
\beq\label{dirichlet.matrix.measure}
C_{d; \bfa}\prod^{n}_{k=1}\det(\bfZ^{(k)})^{a_k-1}\det(\Id - \sum^{n}_{k = 1}\bfZ^{(k)})^{a_{n+1}-1}\prod^{n}_{k= 1}d\bfZ^{(k)},
\eeq
where $\bfa = (a_1, \cdots, a_{n+1})$, $\{a_{i}\}^{n+1}_{i=1}$ are all positive constants, and $d\bfZ^{(k)}$ is the Lebesgue measure on the entries of $\bfZ^{(k)}$, i.e.,  if  $\bfZ^{(k)}$ has complex  entries $Z^{(k)}_{pq}= X^{(k)}_{pq}+ i Y^{(k)}_{pq}$, with $Z^{(k)}_{qp}= \bar Z^{(k)}_{pq}$,  $$d\bfZ^{(k)}= \prod_{1\leq  p\leq q\leq d} dX^{(k)}_{pq}dY^{(k)}_{pq}.$$ This measure is finite exactly when the constants $a_i$ are positive, and the normalizing constant   $ C_{d; \bfa}$ makes $D_{\bfa}$ a probability measure on $\Delta_{n,d}$.

The normalization constant  $ C_{d; \bfa}$  may be explicitly  computed with the help of the matrix gamma function
 $$\Gamma_{d}(a) = \int_{A > 0} e^{-\tr(A)}\det(A)^{a-1}dA = \pi^{\frac{1}{2}d(d-1)}\prod^{d}_{i=1}\Gamma(a+d -i),$$  where  $\{A > 0\}$ denotes  the domain of  $d \times d$ positive-definite, Hermitian matrices. Then,  the normalization constant may be written as
$$C_{d; \bfa} = \frac{\prod^{n+1}_{i=1} \Gamma_{d}(a_{i})}{\Gamma_{d}(\sum^{n+1}_{i=1}a_{i})}.$$ 

It turns out that $\Delta_{n,d}$ is polynomial domain as described in Remark~\ref{rmq.def.polynomial.domain}, with boundary described by the equation
\beqna\label{boundary.eq.mdp}
\det(\bfZ^{(1)}) \cdots \det(\bfZ^{(n)}) \det(\Id- \bfZ^{(1)}-\cdots -\bfZ^{(n)})=0,
\eeqna which is an algebraic equation in the coordinates $(X^{(k)}_{pq}, Y^{(k)}_{pq})$. For convenience, and as described in Section~\ref{sec.polynomial.models}, we will use complex coordinates $(Z^{(k)}_{pq}, \bar Z^{(k)}_{pq},  1\leq k \leq n, 1\leq p <q\leq d)$ (keeping in mind that $Z^{(k)}_{pp}$ which is real),  to describe the various diffusion operators acting on $\Delta_{\bbC, n,d}$. Moreover, instead of using $(Z^{(k)}_{ij}, \bar Z^{(k)}_{ij})$, $1\leq i < j\leq d$ as coordinates,  in our case it is simpler to use $(Z^{(k)}_{ij}$, $i,j=1, \cdots, d)$, due to the fact that $\bar Z^{(k)}_{ij}= Z^{(k)}_{ji}$. As in the scalar case studied in Section~\ref{sec.scalar.simplex}, it turns out that there are many possible operators $\Gamma$ acting on $\Delta_{n,d}$ such that
$(\Delta_{n,d}, \Gamma, D_{a_1, \cdots, a_{n+1} })$ is a polynomial model.

To simplify the notations, we shall always set $\bfZ^{(n+1)}= \Id-\sum_{i=1}^n \bfZ^{(i)}$, with entries $Z^{(n+1)}_{ij}$, $i,j= 1, \cdots, d$.

 According to the boundary equation~\eqref{boundary.eq2}, the carr\'e du champ operator  $\Gamma$  of the matrix Dirichlet process is such that each entry must be a polynomial of degree at most 2  in the variables $(Z_{ij}^{(k)}, \bar Z_{ij}^{(k)})$ and must satisfy
\beqnas\label{beq.dirichlet}
\Gamma(Z^{(p)}_{ij}, \log \det (\bfZ^{(k)})) = L^{p,k},
\eeqnas
for every $p= 1, \cdots, n$, $k= 1, \cdots, n+1$, $i,j = 1, \cdots d$, where  $L^{p,k}$ is an affine  function of the entries of $\{\bfZ^{(1)}, ..., \bfZ^{(n)}\}$. However,  since there are many variables in the diffusion operator $\Gamma$, we are not in the position to describe all the possible solutions for $\Gamma$ as we did in the scalar case.  Therefore, we will restrict ourselves to a simpler form. Namely, we assume that for any $1 \leq p, q \leq n+1$,
$$
\Gamma(Z^{(p)}_{ij}, Z^{(q)}_{kl}) = \sum_{abcdrs}(A^{p, q}_{ij, kl})^{ab, cd}_{r, s}Z^{(r)}_{ab}Z^{(s)}_{cd},
$$ with some constant coefficients $(A^{p, q}_{ij, kl})^{ab, cd}_{r, s}$.
Then the tensor $(\bfA^{p, q}_{ij, kl})$, whose entries are denoted by $\{(A^{p, q}_{ij, kl})^{ab, cd}_{r, s}\}$,  should satisfy the following restrictions:
\bitem
\item{Since $\Gamma$ is symmetric, we have $(A^{p, q}_{ij, kl})_{ab, cd} = (A^{q, p}_{kl, ij})_{ab, cd}$, and the choice of  $(\bfA^{p, q}_{ij, kl})$ should ensure that $\Gamma$ is elliptic on the matrix simplex;
}
\item{The fact that the diffusions live in the matrix simplex gives rise to
$$
\sum_{p=1}^{n+1} \sum_{abcd}(A^{p, q}_{ij, kl})_{ab, cd}Z^{(p)}_{ab}Z^{(q)}_{cd}= 0;
$$
}
\item{The boundary equation \eqref{beq.dirichlet} leads to
\beq\label{eq.GammaLogDelta}
\Gamma(Z^{(p)}_{ij}, \log \det \bfZ^{(q)}) = \sum_{abcd}(Z^{(q)})^{-1}_{lk}(A^{p, q}_{ij, kl})_{ab, cd}Z^{(p)}_{ab}Z^{(q)}_{cd} = L^{p,q},
\eeq
for every $p= 1, \cdots, n$,$q= 1, \cdots, n+1$, where  $L^{p,q}$ is an affine function in the entries of $\{\bfZ^{(1)}, ..., \bfZ^{(n)}\}$, and $(Z^{(q)})^{-1}_{ij}$ are the entries of the inverse matrix $({\bfZ}^{(q)})^{-1}$.

This last equation comes from the diffusion property for the operator $\Gamma$ (equation~\eqref{chain.rule.Gamma}), together with the fact that, for any matrix $\bfZ$ with entries $Z_{ij}$,
$$\partial_{Z_{ij}} \log \det(\bfZ)= Z^{-1}_{ji},$$
where $Z^{-1}_{ij}$ are the entries of the inverse matrix $\bfZ^{-1}$.
}
\eitem

Even under above restrictions, it is still hard  to give any complete description of such tensors  $\{(\bfA^{p, q}_{ij, kl})\}$.  In the following we give two models of matrix Dirichlet process, which appear quite naturally as projections from the Brownian motion on $SU(N)$,  from polar decompositions of  complex matrix Brownian motions, and  from  complex Wishart processes, as we have mentioned before.

\subsection{The first polynomial diffusion model on $\Delta_{n,d}$}\label{sec.model.i}
Our first model is defined in  the following Theorem~\ref{matrix.dirichlet}. As we will see, it appears naturally in some projections of diffusion models on $SU(N)$ and from polar decomposition of complex matrix Brownian motions.

\bthm\label{matrix.dirichlet}
Let the matrix $\bfA = (A_{pq}), 1 \leq p, q \leq n+1$ be a symmetric matrix. Then, consider the  diffusion $\Gamma_{\bfA}$  operator given by
\beq\label{gamma.dirichlet}
\Gamma_{\bfA}(Z^{(p)}_{ij}, Z^{(q)}_{kl}) =  \delta_{pq}\sum^{n + 1}_{s=1} A_{sp}(Z^{(s)}_{il}Z^{(p)}_{kj} +Z^{(s)}_{kj}Z^{(p)}_{il}) - A_{pq}(Z^{(q)}_{il}Z^{(p)}_{kj} +Z^{(p)}_{il}Z^{(q)}_{kj}),
\eeq
for $1\leq p, q \leq n$, $1\leq i,j,k,l\leq d$. Following the notations of Section~\ref{subsec.matrix.Dirichlet.domain}, in this model we have 
$$
(\bfA^{p, q}_{ij, kl})^{ab, cd}_{r, s} = A_{rs}(\delta_{pq}\delta_{pr} - \delta_{pr}\delta_{qs})(\delta_{ak}\delta_{bj}\delta_{ci}\delta_{dl} + \delta_{ai}\delta_{bl}\delta_{ck}\delta_{dj}).
$$
Then $\Gamma_\bfA$ is elliptic if and only if the matrix $A$ has non negative entries and is irreducible.
In this case,  $(\Delta_{n,d}, \Gamma_{\bfA}, D_{\bfa})$ is a polynomial model.  Moreover, in this case,
\beq\label{l.dirichlet}
\LL_{\bfA, \bfa}(Z^{(p)}_{ij})= \sum^{n + 1}_{q }2(a_{p} + d - 1)A_{pq}Z^{(q)}_{ij} - \sum^{n + 1}_{q } 2(a_{q} + d - 1)A_{pq}Z^{(p)}_{ij}.
\eeq
\ethm

We recall that a matrix $A$ with non negative entries is irreducible if and only if for any $p\neq q$, there exists a path $p= p_0, p_1, \cdots, p_k= q$ such that for any $i= 0, \cdots, k-1$, $A_{p_ip_{i+1}} \neq 0$.

\bpf

Recall  that our coordinates are the entries $Z^{(k)}_{ij}$ of the  Hermitian matrices  $\bfZ^{(1)}, \cdots, \bfZ^{(n)}$ and that $\bfZ^{(n + 1)} = \Id - \sum^{n}_{k=1}\bfZ^{(k)}$.


The only requirement here (apart from the ellipticity property  that we willl deal with below) is that equation~\eqref{eq.GammaLogDelta} is satisfied for this particular choice of the tensor $(\bfA^{p,q}_{ij,kl})$. Thus for $1 \leq p, q \leq n$, we have
$$
\Gamma_{\bfA}(\log \det(\bfZ^{(q)}), Z^{(p)}_{ij}) = \delta_{pq}\sum^{n+1}_s2A_{sp}Z^{(s)}_{ij} - 2A_{pq}Z^{(p)}_{ij} ,$$
while
$$
\Gamma_{\bfA}(\log \det(\bfZ^{(n+1)}), Z^{(p)}_{ij})
=  -2A_{(n+1)p}Z^{(p)}_{ij},
$$
which shows that the boundary equation is satisfied for this model.

We now  prove that on the matrix simplex  $\Delta_{n,d}$ , $\Gamma_{\bfA}$ given by \eqref{gamma.dirichlet} is elliptic if and only if for $p, q = 1,...., n+1$, $A_{pq} > 0$.


For  fixed $(p, q)$,  we consider $\Gamma_{\bfA}(Z^{(p)}_{ij}, Z^{(q)}_{kl})$ as the $(ij, kl)$ element in a $d^{2} \times d^{2}$ matrix;  and $(\Gamma_{\bfA}(Z^{(p)}_{ij}, Z^{(q)}_{kl}))$ is the $(p, q)$ element in a  $n\times n$ matrix of $d^2\times d^2$ matrices. Notice that we may write 
$$
(\Gamma_{\bfA}) =  \sum_{1\leq p<q\leq n}A_{pq}\Gamma^{(pq)} + \bfA_{n+1}\Gamma^{(n+1)},
$$
where 

\benum 

\item $(\Gamma^{pq})$ is a  $nd^{2} \times nd^{2}$ block matrix with
\beqnas
(\Gamma^{(pq)})_{pq, (ij, kl)} &=& - (Z^{(q)}_{il}Z^{(p)}_{kj} +Z^{(p)}_{il}Z^{(q)}_{kj}),\\
(\Gamma^{(pq)})_{qp, (ij, kl)} &=& - (Z^{(q)}_{il}Z^{(p)}_{kj} +Z^{(p)}_{il}Z^{(q)}_{kj}),\\
(\Gamma^{(pq)})_{pp, (ij, kl)} &=& Z^{(q)}_{il}Z^{(p)}_{kj} +Z^{(p)}_{il}Z^{(q)}_{kj},\\
(\Gamma^{(pq)})_{qq, (ij, kl)} &=& Z^{(q)}_{il}Z^{(p)}_{kj} +Z^{(p)}_{il}Z^{(q)}_{kj},
\eeqnas
and  all other entries are $0$; 

\item$\Gamma^{(n+1)}$ is a diagonal block matrix with $\Gamma^{(n+1)}_{pp, (ij, kl)} = Z^{(n+1)}_{il}Z^{(p)}_{kj} +Z^{(p)}_{il}Z^{(n+1)}_{kj}$ and all other entries are $0$; 

\item $\bfA_{n+1}$ is a diagonal block  matrix of size $nd^{2} \times nd^{2}$ satisfying $(\bfA_{n+1})_{pp} = A_{p(n+1)}\Id_{d^{2} \times d^{2}}$.

\eenum 

We first observe that each operator associated with $\Gamma^{(pq)}$ is non negative. To see this, we have to check that for any sequence $(\Lambda)= (\Lambda^{(1)}, \cdots, \Lambda^{(n)})$ of Hermitian  matrices   with entries $\lambda^{(k)}_{ij}$,
$$\sum_{rs, ij, kl} \Gamma^{(pq)}_{rs, (ij, \bar{kl})}\lambda^{(r)}_{ij} \bar \lambda^{(s)}_{kl}  \geq 0.$$

In fact,
\beqnas
& &\sum_{rs, ij, kl} \Gamma^{(pq)}_{rs, (ij, \bar{kl})}\lambda^{(r)}_{ij} \bar \lambda^{(s)}_{kl} \\
&=& \Lambda^{(p)}\bar{\Lambda}^{(p)}\Gamma^{(pq)}_{pp} + \Lambda^{(q)}\bar{\Lambda}^{(q)}\Gamma^{(pq)}_{qq} + \Lambda^{(p)}\bar{\Lambda}^{(q)}\Gamma^{(pq)}_{pq} + \Lambda^{(q)}\bar{\Lambda}^{(p)}\Gamma^{(pq)}_{qp}\\
&=& \Lambda^{(p)}\bar{\Lambda}^{(p)}\Gamma^{(pq)}_{pp}  + \Lambda^{(q)}\bar{\Lambda}^{(q)}\Gamma^{(pq)}_{pp} - \Lambda^{(p)}\bar{\Lambda}^{(q)}\Gamma^{(pq)}_{pp}  - \Lambda^{(q)}\bar{\Lambda}^{(p)}\Gamma^{(pq)}_{pp} \\
&=& (\Lambda^{(p)} - \Lambda^{(q)})\Gamma^{(pq)}_{pp}(\bar{\Lambda}^{(p)} - \bar{\Lambda}^{(q)})\\
&=& \sum_{ij, kl} (\lambda^{(p)} - \lambda^{(q)})_{ij}(\bar{\lambda}^{(p)} - \bar{\lambda}^{(q)})_{kl}(Z^{(q)}_{ik}Z^{(p)}_{lj} + Z^{(p)}_{ik}Z^{(q)}_{lj})\\
&=& \tr(\bfZ^{(q)}(\bar{\Lambda}^{(p)} - \bar{\Lambda}^{(q)})\bfZ^{(p)}(\Lambda^{(p)} - \Lambda^{(q)})^{t}) \\
& & \ \ \ \ + \tr(\bfZ^{(p)}(\bar{\Lambda}^{(p)} - \bar{\Lambda}^{(q)})\bfZ^{(q)}(\Lambda^{(p)} - \Lambda^{(q)})^{t}).
\eeqnas

For any $p$, $\bfZ^{(p)}$ is Hermitian and non negative-definite, so are $(\bar{\Lambda}^{(p)} - \bar{\Lambda}^{(q)})\bfZ^{(p)}(\Lambda^{(p)} - \Lambda^{(q)})^{t}$, $(\Lambda^{(p)} - \Lambda^{(q)})^{t}\bfZ^{(p)}(\bar{\Lambda}^{(p)} - \bar{\Lambda}^{(q)})$, then we have
\beqnas
\tr(\bfZ^{(q)}(\bar{\Lambda}^{(p)} - \bar{\Lambda}^{(q)})\bfZ^{(p)}(\Lambda^{(p)} - \Lambda^{(q)})^{t}) &\geq& 0, \\
\tr(\bfZ^{(p)}(\bar{\Lambda}^{(p)} - \bar{\Lambda}^{(q)})\bfZ^{(q)}(\Lambda^{(p)} - \Lambda^{(q)})^{t}) &\geq& 0.
\eeqnas
therefore $\Gamma^{(pq)}$ are all non negative-definite matrices. 

Similarly, we have
\beqnas
& &\sum_{p, ij, kl} A_{p(n+1)}\Gamma^{(n+1)}_{pp, (ij, \bar{kl})}\lambda^{(p)}_{ij} \bar \lambda^{(p)}_{kl} \\
&=& \sum_{p, ij, kl} A_{p(n+1)}(Z^{(n+1)}_{ik}Z^{(p)}_{lj} +Z^{(p)}_{ik}Z^{(n+1)}_{lj})\lambda^{(p)}_{ij} \bar \lambda^{(p)}_{kl} \\
&=& \sum_{p} A_{p(n+1)}(\tr(\bfZ^{(n+1)}\bar \Lambda^{(p)}\bfZ^{(p)}(\Lambda^{(p)})^{t}) + \tr(\bfZ^{(n+1)}(\Lambda^{(p)})^{t}\bfZ^{(p)}\bar \Lambda^{(p)})),
\eeqnas
Thus we know $\Gamma^{(n+1)} \geq 0$.

We then prove the following lemma,
 \blem\label{positive-definite}
 Let $A$, $B$ be $d \times d$ Hermitian positive definite matrices. Then for a given $d \times d$ matrix $U$, if $\tr(AUBU^*)=0$, then $U=0$.
 \elem

\bpf
Suppose $A$ has the spectral decomposition $A = P^{*}DP$, where $P$ is unitary and $D = \diag\{\lambda_1, \cdots, \lambda_d\}$ with all $\lambda_{i}$ positive. Then
\beqnas
\tr(AUBU^*) = \tr(P^{*}DPUBU^*) = \tr(DPUBU^*P^{*}).
\eeqnas 
Notice that since $B$ is positive definite, if $PU \neq 0$, then $PUBU^*P^{*}$ is also positive definite, whose elements on the diagonal are all positive, implying $\tr(DPUBU^*P^{*}) > 0$. 
 Therefore, $\tr(AUBU^*) =0$ holds only when $PU = 0$. Since $P$ is unitary, this happens only when $U = 0$.
\epf

Recall the boundary equation of $\Delta_{n,d}$ \eqref{boundary.eq.mdp}, then we know all $\bfZ^{(i)}$ are positive definite inside $\Delta_{n,d}$ .

Now we claim that $\Gamma_{\bfA}$ is elliptic inside $\Delta_{n,d}$ if and only if $\bfA$ has non negative entries and is irreducible. 

Previous computations show that $\Gamma_{\bfA} \geq 0$.  If there exists $\Lambda_{0}$ such that $\Lambda_{0} \Gamma_{\bfA} \Lambda_{0}^{*} = 0$, we have
\beqna\label{ellipticity}
 & &\sum^{n}_{p < q}A_{pq}(\tr(\bfZ^{(q)}(\bar{\Lambda}_{0}^{(p)} - \bar{\Lambda}_{0}^{(q)})\bfZ^{(p)}(\Lambda^{(p)}_{0} - \Lambda_{0}^{(q)})^{t}) \\
\nonumber& & \ \ \ \ + \tr(\bfZ^{(p)}(\bar{\Lambda}_{0}^{(p)} - \bar{\Lambda}_{0}^{(q)})\bfZ^{(q)}(\Lambda_{0}^{(p)} - \Lambda^{(q)})_{0}^{t})) \\
\nonumber & & \ \ \ \ + \sum^{n}_{p}A_{p(n+1)}(\tr(\bfZ^{(n+1)}\bar \Lambda^{(p)}\bfZ^{(p)}(\Lambda^{(p)})^{t}) \\
\nonumber& & \ \ \ \ + \tr(\bfZ^{(n+1)}(\Lambda^{(p)})^{t}\bfZ^{(p)}\bar \Lambda^{(p)})) = 0.
\eeqna
Since  all the  terms in the previous sum are non-negative, they  all vanish. By the fact that $\bfA$ is irreducible, we know that for any $1 \leq p \leq n$, there exists a path connecting $p$ and $n+1$, noted by $p, p_{1}, p_{2}, \cdots, p_{r}, (n+1)$, which leads to $A_{pp_{1}}, A_{p_{1}p_{2}}, \cdots, A_{p_{r}(n+1)} > 0$. Then by Lemma \ref{positive-definite}, we have
$$
\Lambda^{(p)}_{0} = \Lambda^{(p_{1})}_{0} = \cdots = \Lambda^{(p_{r})}_{0} = 0,
$$  
i.e. $\Lambda_{0} = 0$, which implies that $\Gamma_{\bfA}$ is elliptic. 

Conversely, when $\Gamma_{\bfA}$ is elliptic, then if $\Lambda_{0} \Gamma_{\bfA} \Lambda_{0}^{*} = 0$, we must have $\Lambda_{0} = 0$. First we prove that $\bfA$ has non negative entries. For given $1 \leq p \leq (n+1)$, assume that some $A_{pq}<0$ for some $q\neq p$. Then, choose the sequence  $\Lambda_{1} = (\Lambda^{(1)}_{1}, \cdots , \Lambda_{1}^{(n)})$ to be such that $\Lambda^{(p)} = ((\bfZ^{(p)})^{-\frac{1}{2}})^{t}$ and all others are $0$, then from $\Lambda_{1}\Gamma_{\bfA}\Lambda_{1}^{*} > 0$,
we conclude that for $(\bfZ^{(1)}, \cdots, \bfZ^{(n)}) \in \Delta_{n,d}$, 
\beqnas
\sum^{n+1}_{r \neq p} 2A_{pr}\tr(\bfZ^{(r)}) > 0.
\eeqnas

Now, choose $Z^{(q)}= (1-\epsilon)\Id$ and $Z^{(r)}= \frac{\epsilon}{n} \Id$ for $r\neq q$, $0 < \epsilon < 1$. 
Then when $\epsilon < \frac{|A_{pq}|}{\sum^{n+1}_{s \neq p}|A_{ps}|}$  we obtain a contradiction.

 Now suppose $\bfA$ is not irreducible, then it has at least two strongly connected components $\cA_{1}$, $\cA_{2}$, such that $A_{pq} = 0$ for $p \in \cA_{1}$, $q \in \cA_{2}$.  Suppose $(n+1) \in \cA_{1}$, then for all $p \in \cA_{1}$, choose $\Lambda^{(p)}_{0} = 0$, while for $q \in \cA_{2}$, $\Lambda^{(q)}_{0} \neq 0$, such that we have $\Lambda_{0} \neq 0$ which satisfies $\Lambda_{0} \Gamma_{\bfA} \Lambda_{0}^{*} = 0$. Thus there is a contradiction, so we must have $\bfA$ to be irreducible.

Finally, by a direct computation, we have
\beqnas
\LL_{\bfA, \bfa}(Z^{(p)}_{ij})
&=& \sum^{n+1}_{q=1}(a_{q} - 1)\Gamma(\log \det(\bfZ^{(q)}), Z^{(p)}_{ij}) + \sum^{n}_{q=1}\sum_{kl}\partial_{Z^{(q)}_{kl}}\Gamma(Z^{(p)}_{ij}, Z^{(q)}_{kl})\\
&=& \sum^{n+1}_{q }2(a_{p} + d - 1)A_{pq}Z^{(q)}_{ij} - \sum^{n+1}_{q } 2(a_{q} + d - 1)A_{pq}Z^{(p)}_{ij}.
\eeqnas
\epf

In the sequel, we give two natural constructions of matrix Dirichlet process, the first one  from the  Brownian motion on $SU(N)$, the other from polar decomposition of  Brownian complex matrices.

\subsubsection{The construction from $SU(N)$}\label{construction2}
In this section, we show that the matrix Dirichlet processes can be realized by the projections from Brownian motion on $SU(N)$. The construction relies on the matrix-extracting procedure,  extending the construction of matrix Jacobi processes described in~\cite{Doumerc}. In  the special case where all the parameters $A_{ij}$ are equal to $d$, the associated generator may be considered as an image of  the Casimir operator on $SU(N)$, whenever the coefficients $a_i$ in the measure are integers. Moreover, for the general case, similar to  the construction  in the scalar case in Proposition~\ref{prop.general.scalar.case}, we provide a construction from more general Brownian motions on $SU(N)$, where the generator is no longer the Casimir operator on $SU(N)$.

The so-called "matrix-extracting procedure" is the following: consider a matrix $u$ on $SU(N)$, then take the first $d$ lines, and split the set of all $N$ columns into $(n + 1)$ disjoint sets $I_1, ..., I_{n +1}$. For $1 \leq i \leq (n+1)$, define $d_i = |I_i|$ such that $N = d_1+\cdots + d_{n+1}$. Then we get $(n+1)$ extracted matrices $\{\bfW^{(i)}\}$, respectively of size $d \times d_1$, $d \times d_{2}$,..., $d \times d_{n +1}$.  The matrix Jacobi process is obtained by considering the Hermitian matrix $\bfW^{(1)}\bfW^{(1)}{}^*$. Here, we extend  this procedure,  defining $\bfZ^{(i)} = \bfW^{(i)}(\bfW^{(i)})^{*}$ for $i =1, ..., n+1$. Then, $\sum_{i=1}^{n+1} \bfZ^{(i)}= \Id$, and the process $(\bfZ^{(1)}, \cdots, \bfZ^{(n)})$ lives in the complex matrix Dirichlet simplex $\Delta_{d, n}$.

The compact Lie group $SU(N)$ is semi-simple compact. There is, up to a scaling constant, a unique elliptic diffusion operator on it  which commutes both with the right and the left multiplication. This operator is called the Casimir operator, see \cite{BOZ}, \cite{zribi} for more details. The Brownian motion on $SU(N)$ is the diffusion process which has the Casimir operator as its generator. It may be described by the vector fields  $\cV_{R_{ij}}$, $\cV_{S_{ij}}$ and $\cV_{D_{ij}}$ which are given on the entries $\{u_{ij}\}$ of  $\bfu \in SU(N)$ matrix as
\beqna
\label{casimir.r}\cV_{R_{ij}} &=& \sum_{k}\big(u_{kj}\partial_{u_{ki}} - u_{ki}\partial_{u_{kj}} + \bar{u}_{kj}\partial_{\bar{u}_{ki}} - \bar{u}_{ki}\partial_{\bar{u}_{kj}}\big),\\
\label{casimir.s}\cV_{S_{ij}} &=& i\sum_{k}\big( u_{kj}\partial_{u_{ki}} + u_{ki}\partial_{u_{kj}}- \bar{u}_{kj}\partial_{\bar{u}_{ki}} - \bar{u}_{ki}\partial_{\bar{u}_{kj}} \big),\\
\label{casimir.d}\cV_{D_{ij}} &=& i\sum_{k}\big(u_{ki}\partial_{u_{ki}} - u_{kj}\partial_{u_{kj}} - \bar{u}_{ki}\partial_{\bar{z}_{ki}} + \bar{u}_{kj}\partial_{\bar{u}_{kj}}\big).
\eeqna
Then for the Casimir operator $\Delta_{SU(N)}$  on $SU(N)$, we have

\beqnas\label{casimir.opt}
\Delta_{SU(N)} = \frac{1}{4N}\sum_{i < j}(\cV^{2}_{R_{ij}} + \cV^{2}_{S_{ij}} + \frac{2}{N}\cV^{2}_{D_{ij}}),
\eeqnas
and
\beqnas
\Gamma_{SU(N)}(f, g) &=& \frac{1}{4N}(\sum_{i < j}\cV_{R_{ij}}(f)\cV_{R_{ij}}(g) + \cV_{S_{ij}}(f)\cV_{S_{ij}}(g) + \frac{2}{N}\cV_{D_{ij}}(f)\cV_{D_{ij}}(g)).
\eeqnas

From this, a simple computation provides
\beqnas
\Gamma_{SU(N)}(u_{ij}, u_{kl}) &=& -\frac{1}{2N}u_{il}u_{kj} + \frac{1}{2N^{2}}u_{ij}u_{kl},\\
\Gamma_{SU(N)}(u_{ij}, \bar{u}_{kl}) &=& \frac{1}{2N}\delta_{ik}\delta_{jl} - \frac{1}{2N^2}u_{ij}\bar{u}_{kl},\\
\Delta_{SU(N)}(u_{ij}) &=& -\frac{N^{2}-1}{2N^{2}}u_{ij}, \ \ \ \ \Delta_{SU(N)}(\bar{u}_{ij}) = -\frac{N^{2}-1}{2N^{2}}\bar{u}_{ij},\\
\eeqnas
and these formulas describe entirely the Brownian motion on $SU(N)$.

A simple application of the diffusion property (equations~\eqref{chain.rule.L} and \eqref{chain.rule.Gamma}) yields, for $1 \leq p, q \leq n$, for the entries $Z^{(p)}_{ij}$ of the matrices $\bfZ^{(p)}$,
\beqna
 & &\label{gamma.ug1}\Gamma_{SU(N)}(Z^{(p)}_{ij}, Z^{(q)}_{kl}) = \frac{1}{2N}\delta_{pq}(\delta_{il}Z^{(p)}_{kj}  + \delta_{kj}Z^{(p)}_{il}) - \frac{1}{2N}(Z^{(p)}_{kj}Z^{(q)}_{il} + Z^{(p)}_{il}Z^{(q)}_{kj}),  \\
 & &\label{l.ug2}\Delta_{SU(N)}(Z^{(p)}_{ij}) = - Z^{(p)}_{ij} + \frac{1}{N}d_{p}\delta_{ij}.
\eeqna
Comparing \eqref{gamma.ug1}, \eqref{l.ug2} with \eqref{gamma.dirichlet} and \eqref{l.dirichlet}, we obtain a matrix Dirichlet operator described in Theorem~\ref{matrix.dirichlet}, with
$$
A_{pq}= \frac{1}{2N}, \ \ \ a_i = d_i - d + 1,
$$
for any $1 \leq p, q \leq  n$ and $1 \leq i \leq (n + 1)$. Therefore, the density of the reversible measure is
\beq
C\prod^{n}_{p = 1}\det(\bfZ^{(p)})^{d_{p} - d}\det(\Id - \sum^{n}_{p=1}\bfZ^{(p)})^{d_{n+1} - d}.
\eeq
For  this measure  to be finite,  we need $d_i > d -1$ for all $1 \leq i \leq (n+1)$, i.e., $d_i \geq d$ since these parameters are integers. It is worth to observe that this restriction is necessary for the matrices $\bfZ^{(i)}$ to be non degenerate. If it is not satisfied, the matrices $Z^{(i)}$ live of an algebraic manifold and their law may not have any density with respect to the Lebesgue measure.

To summarize, we have

\bprop The image of the Brownian motion on $SU(N)$ under the matrix-extracting procedure is a diffusion process on the complex matrix simplex $\Delta_{d,n}$ with carré du champ  $\Gamma_{\bfA}$ and  reversible measure $D_{a_1, \cdots, a_{n+1}}$ where $A_{pq}=  \frac{1}{2N}, p, q= 1, \cdots, n+1, p\neq q$, 
$a_i = d_i -d + 1$.
\eprop

As a corollary, we get
\bcor Whenever $d_i \geq d$, $i= 1, \cdots, n+1$, the image of the Haar measure on $SU(N)$ through the matrix-extracting procedure is the matrix Dirichlet measure $D_{d_1-d, \cdots, d_{n+1}-d}$.
\ecor

For the general case where the parameters $A_{ij}$ are not equal, we may follow Proposition~\ref{prop.general.scalar.case}.
For $p \neq q$, we define the following $\LL^{(pq)}$ acting on the matrix simplex \\ $\{\bfZ^{(i)}, \sum^{n}_{i=1} \bfZ^{(i)} \leq \Id\}$,
$$\LL^{ (pq)}= \sum_{i\in I_p, j\in I_q } \cV_{R_{ij}}^2+ \cV_{S_{ij}}^2+ \frac{2}{N} \cV_{D_{ij}}^2,$$
with its corresponding carré du champ operator $\Gamma^{(p q)}$,
\beqnas
\Gamma^{(p q)}(f, g) = \sum_{i \in I_p,  j \in I_q}\cV_{R_{ij}}(f)\cV_{R_{ij}}(g) + \cV_{S_{ij}}(f)\cV_{S_{ij}}(g) + \frac{2}{N}\cV_{D_{ij}}(f)\cV_{D_{ij}}(g).
\eeqnas

\blem\label{xeq0}
For the entries $u_{ij}$ of an $SU(N)$ matrix, and denoting by  $Z^{(p)}_{ij}$ the entries of the extracted matrix $\bfZ^{(p)}$, we have
\beqna
 \label{xeq1}\cV_{R_{ij}}Z^{(p)}_{ab} &=& \un_{i \in I_p}(u_{aj}\bar{u}_{bi} + u_{ai}\bar{u}_{bj})- \un_{j \in I_p}(u_{ai}\bar{u}_{bj}  + u_{aj}\bar{u}_{bi}),\\
 \label{xeq2}\cV_{S_{ij}}Z^{(p)}_{ab}  &=& \sqrt{-1}\big(\un_{i \in I_p}(u_{aj}\bar{u}_{bi}  - u_{ai}\bar{u}_{bj}) + \un_{j \in I_p}(u_{ai}\bar{u}_{bj} - u_{aj}\bar{u}_{bi})\big),\\
 \label{xeq3}\cV_{D_{ij}}Z^{(p)}_{ab} &=& 0,
\eeqna

and
\beqna
\label{xeq8}\cV^{2}_{R_{ij}}Z^{(p)}_{ab} &=& \un_{i  \in I_p}(2u_{aj}\bar{u}_{bj} - 2u_{ai}\bar{u}_{bi}) - \un_{j \in I_p}(2u_{aj}\bar{u}_{bj} - 2u_{ai}\bar{u}_{bi}),\\
 \label{xeq9}\cV^{2}_{S_{ij}}Z^{(p)}_{ab} &=& \un_{i  \in I_p}(2u_{aj}\bar{u}_{bj} - 2u_{ai}\bar{u}_{bi}) - \un_{j \in I_p}(2u_{aj}\bar{u}_{bj} - 2u_{ai}\bar{u}_{bi}).
\eeqna

Then,
\beqna
\label{xeq4}\Gamma^{(p q)}(Z^{(p)}_{ab}, Z^{(q)}_{cd}) &=& - 2Z^{(p)}_{cb}Z^{(q)}_{ad} - 2Z^{(p)}_{ad}Z^{(q)}_{cb} , \\
\label{xeq5}\Gamma^{(p q)}(Z^{(p)}_{ab}, Z^{(p)}_{cd})  &=&  2Z^{(p)}_{cb}Z^{(q)}_{ad} + 2Z^{(p)}_{ad}Z^{(q)}_{cb},  \\
\label{xeq6}\Gamma^{(p q)}(Z^{(q)}_{ab}, Z^{(q)}_{cd}) &=& 2Z^{(p)}_{cb}Z^{(q)}_{ad} + 2Z^{(p)}_{ad}Z^{(q)}_{cb}.
\eeqna

For a pair $(r, s) \neq (p, q)$, we have
\beqna
\label{xeq7}\Gamma^{(pq)}(Z^{(r)}_{ab}, Z^{(s)}_{cd}) = 0.
\eeqna

Moreover,
\beqna
 \label{xeq10}\LL^{(p q)}(Z^{(r)}_{ij}) &=& \un_{r=p}4(d_{p}Z^{(q)}_{ij} - d_{q}Z^{(p)}_{ij}) - \un_{r=q}4(d_{p}Z^{(q)}_{ij} - d_{q}Z^{(p)}_{ij}).
\eeqna
\elem
 
\brmq 
 The action of $\cV_{D_{ij}}$ vanishes  on the variables $Z^{(p)}_{ab}$ such that indeed we could replace $\LL^{(pq)}$ by
$$\LL^{ (pq)}= \sum_{i\in I_p, j\in I_q } \cV_{R_{ij}}^2+ \cV_{S_{ij}}^2,$$ and the dimension $N$ disappears  from the definition.
\ermq

\bpf
Recall \eqref{casimir.r}. Then letting $\cV_{R_{ij}}$ act on $\bfZ^{p} = \bfW^{p}(\bfW^{p})^{*}$ and writing $u_{ij}$ for the entries of $\bfW^{(p)}$,  which are indeed entries of an $SU(N)$ matrix, we have
\beqnas
\cV_{R_{ij}}Z^{(p)}_{ab}
&=&  \sum_{r \in I_p}\cV_{R_{ij}}(u_{ar}\bar{u}_{br})\\
&=& \un_{i \in I_p}(u_{aj}\bar{u}_{bi} + u_{ai}\bar{u}_{bj})- \un_{j \in I_p}(u_{ai}\bar{u}_{bj}  + u_{aj}\bar{u}_{bi}).
\eeqnas

Similarly we can prove \eqref{xeq2}, \eqref{xeq3}.

By \eqref{xeq1}, \eqref{xeq2}, we obtain
\beqnas
\cV_{R_{ij}}(u_{aj}\bar{u}_{bi} + u_{ai}\bar{u}_{bj})
&=& 2u_{aj}\bar{u}_{bj} - 2u_{ai}\bar{u}_{bi},\\
\cV_{S_{ij}}(u_{aj}\bar{u}_{bi} - u_{ai}\bar{u}_{bj})
&=& -\sqrt{-1}(2u_{aj}\bar{u}_{bj} - 2u_{ai}\bar{u}_{bi}),
\eeqnas
then \eqref{xeq7}, \eqref{xeq8} follow.

By \eqref{xeq1}, \eqref{xeq2} and \eqref{xeq3} we have
\beqnas
& &\Gamma^{(pq)}(Z^{(p)}_{ab}, Z^{(q)}_{cd}) \\
&=& \sum_{i \in I_p, j \in I_q}  \cV_{R_{ij}}(Z^{(p)}_{ab})\cV_{R_{ij}}(Z^{(q)}_{cd}) + \cV_{S_{ij}}(Z^{(p)}_{ab})\cV_{S_{ij}}(Z^{(q)}_{cd})\\
&=& - 2Z^{(p)}_{cb}Z^{(q)}_{ad} - 2Z^{(p)}_{ad}Z^{(q)}_{cb} ,
\eeqnas
which proves \eqref{xeq4}.

In the same way, we can deduce \eqref{xeq5}, \eqref{xeq6} and \eqref{xeq7}.

\eqref{xeq10} is proved by \eqref{xeq8}, \eqref{xeq9} and
\beqnas
\LL^{(p q)}(Z^{(r)}_{ij}) &=& \sum_{k \in I_p, l \in I_{q}}\cV^{2}_{R_{kl}}Z^{(r)}_{ij} + \cV^{2}_{S_{kl}}Z^{(r)}_{ij}\\
&=& \un_{r=p}4(d_{p}Z^{(q)}_{ij} - d_{q}Z^{(p)}_{ij}) - \un_{r=q}4(d_{p}Z^{(q)}_{ij} - d_{q}Z^{(p)}_{ij}).
\eeqnas
\epf

Now by Lemma \ref{xeq0} we may derive the following conclusion. 
\bprop
 $(\mathbf{Z}^{(1)}, \cdots, \mathbf{Z}^{(n)})$ is a closed system for any $\LL^{(pq)}$ and the image of
 $$
 \frac{1}{2}\sum^{n+1}_{p < q} A_{pq}\LL^{(pq)}
 $$
 is the operator $\LL_{\bfA, \bfa}$ in Theorem \ref{matrix.dirichlet} with $a_i = d_i -d + 1$.
\eprop

\subsubsection{The construction from polar decomposition}\label{construction1}
In this section, we consider the spectral decomposition of a complex Brownian matrix $\bfm$. Any complex matrix $\bfm$ may be written as  $\bfm = \bfV \bfN$, where $\bfV$ is unitary and $\bfN$ is non negative  Hermitian,  and uniquely determined  as $\sqrt{\bfm^{*}\bfm}$.

A complex Brownian matrix is a diffusion process whose generator may be described on the entries $m_{ij}$ of the $d \times d$ complex matrix $\bfm$ by
$$
\Gamma(m_{ij}, m_{kl}) = 0, \ \ \Gamma(m_{ij}, \bar{m}_{kl}) = 2\delta_{ik}\delta_{jl}, \ \ \ \LL(m_{ij}) = 0.
$$
This just describes the fact that the entries $m_{ij}$ are independent complex Brownian motions. The polar decomposition of a complex Brownian matrix will be discussed in Appendix~\ref{polar}.

Now let $\bfH = \bfm^{*}\bfm$ and suppose $\bfH$ has the spectral decomposition $\bfH = \bfU \bf\Sigma \bfU^{*}$, where $\bfU$ is unitary and $\bf\Sigma$ is  diagonal, denoted by $\bf\Sigma = \diag\{X_1, ..., X_d\}$. Since $\bfH$ is positive definite, we may write $\bf \Sigma = \bfD^{2}$, where $\bfD = \diag \{x_1, \cdots, x_d\}$,  $x_i \geq 0$ and $X_{i} = x^{2}_{i}$ for $1 \leq i \leq d$.
Write $\bfU = (U^{(1)}, \cdots, U^{(d)})$, where $U^{(i)}, 1\leq i \leq d$ is a $d-$column vector. Define $\bfZ^{(k)} = U^{(k)}(U^{(k)})^{*}$, such that for fixed $k$, $\bfZ^{(k)}$ is a Hermitian matrix, and also
\beqnas
\sum^{d}_{k =1}\bfZ^{(k)} = \Id.
\eeqnas

\bprop
The diffusion operators of $(\bfD, \bfZ^{(1)}, \cdots, \bfZ^{(d-1)})$ are given by 
\beqnas
\sum^{d}_{i=1} \big(\partial^{2}_{x_i} + (\frac{1}{x_i} + 4x_i\sum_{j \neq i}\frac{1}{x^{2}_{i} - x^{2}_{j}})\partial_{x_i} \big) + \LL_{\bfA},
\eeqnas
where $\LL_{\bfA}$ is described in Theorem \ref{matrix.dirichlet} with
\beqnas
A_{pq} =  2\frac{x^{2}_{p} + x^{2}_{q}}{(x^{2}_{p} - x^{2}_{q})^2}, \ \ a_{i} = -d,
\eeqnas
for $1 \leq p, q, i \leq d$.

Then the reversible measure is integrable only when $d=1$, which is indeed the Lebesgue measure on the simplex, corresponding to the complex scalar Dirichlet process.  
\eprop

Indeed, the matrices $\bfZ^{(i)}$ have rank one, and therefore their law may not have a density with respect to the Lebesgue measure. This model for matrix Dirichlet processes is in fact degenerate and lives on the boundary of the domain  $\Delta_{n,d}$.

\bpf
The computations of  Appendix~\ref{polar} provide  
\beqnas
\label{ceq1}\Gamma(Z^{(p)}_{ij}, Z^{(q)}_{kl}) &=& -r_{pq}(Z^{(q)}_{il}Z^{(p)}_{kj} +Z^{(p)}_{il}Z^{(q)}_{kj}) + \delta_{pq}\sum^{d}_{s = 1}r_{sp}(Z^{(s)}_{il}Z^{(p)}_{kj} +Z^{(s)}_{kj}Z^{(p)}_{il}), \\
\label{ceq2}\LL(Z^{(p)}_{ij}) &=& 2\sum^{d}_{1 \leq q  \neq p}r_{pq}(Z^{(p)}_{ij} - Z^{(q)}_{ij}),
\eeqnas
where $r_{pq} =  2\frac{X_{p} + X_{q}}{(X_{p} - X_{q})^2}$, which the expected result. 
\epf

In the general case, when $d>1$,  we may  consider the diffusion of $\{v^{(k)} = Z^{(k)}_{11}, 1 \leq k \leq d-1\}$, which is not degenerate and is indeed the complex scalar Dirichlet process.

\brmq
In fact, since $\bfN = \sqrt{\bfH}$, we may write $\bfN = \bfU \bfD \bfU^{*}$, where $\bfD = \sqrt{\bf\Sigma}$. This leads to $\bfm = \bfW \bfD \bfU^{*}$, where $\bfW = \bfV \bfU$. Now we may define $\{\bfY^{(k)} = W^{(k)}(W^{(k)})^{*},  1 \leq k \leq d-1 \}$, where $W^{k}$ is the $k$-th column vector in $\bfW$. Then following Appendix \ref{polar}, we have for  the entries $Y^{(p)}_{ij}$ of $\bfY^{(p)}$
\beqnas
\Gamma(Y^{(p)}_{ij}, Y^{(q)}_{kl}) &=& -\omega_{pq}(Y^{(q)}_{il}Y^{(p)}_{kj} +Y^{(p)}_{il}Y^{(q)}_{kj}) + \delta_{pq}\sum^{d}_{s = 1}\omega_{sp}(Y^{(s)}_{il}Y^{(p)}_{kj} +Y^{(s)}_{kj}Y^{(p)}_{il}), \\
\LL(Y^{(p)}_{ij}) &=& 2\sum^{d}_{1 \leq q  \neq p}\omega_{pq}(Y^{(q)}_{ij} - Y^{(p)}_{ij}),
\eeqnas
where $\omega_{pq} = r_{pq}$ for $p \neq q$, showing that it is exactly the same situation as $\{\bfZ^{(k)}, 1 \leq k \leq d-1\}$.

Although the generators of $\bfW$ are different from those of $\bfU$, as we will see in Appendix A, they have no influence on the diffusion operators of $\bfY^{(k)}$ and $\bfZ^{(k)}$, because the difference lies in $\omega_{pp}$ for $p =1, ..., d$, which make no contribution in the above formulas.

It is also worth to point out that if we consider the left polar decomposition (the previous one is known as the right polar decomposition), i.e., $\bfm = \bfN' \bfV$, where $\bfV$ is the same unitary matrix and $\bfN' = \sqrt{\bfm \bfm^{*}}$,  then $\bfW$ can be viewed as the unitary part in the spectral decomposition of $\bfH' = \bfm \bfm^{*} = \bfW \bf\Sigma \bfW^{*}$.
\ermq

\subsection{The second polynomial diffusion model on $\Delta_{n,d}$}\label{sec.model.ii}

The second model is given by the following Theorem~\ref{gamma.dirichlet2}. As in the scalar case (Section \ref{sec.scalar.simplex}), it can be naturally derived from the Ornstein-Uhlenbeck process on  complex matrices. 

\bthm\label{model2}
Let $\bfA$ be a $d \times d$ positive-definite Hermitian matrix and  $\bfB$ be a $d^{2} \times d^{2}$ positive-definite Hermitian matrix. Consider the diffusion $\Gamma_{\bfA, \bfB}$ operator given by 
\beqna\label{gamma.dirichlet2}
& &\Gamma_{\bfA, \bfB}(Z^{(p)}_{ij}, Z^{(q)}_{kl}) = \delta_{pq}(A_{il}Z^{(p)}_{kj} + A_{kj}Z^{(p)}_{il}) - A_{kj}(Z^{(p)}Z^{(q)})_{il} - A_{il}(Z^{(q)}Z^{(p)})_{kj}\\
\nonumber  & & \ \ \ \ + \sum_{ab}\big( B_{ia, lb}Z^{(p)}_{aj}Z^{(q)}_{kb} +  B_{aj, bk}Z^{(p)}_{ia}Z^{(q)}_{bl} - B_{aj, lb}Z^{(p)}_{ia}Z^{(q)}_{kb} - B_{ia, bk}Z^{(p)}_{aj}Z^{(q)}_{bl} \big) ,  
\eeqna
for $1\leq p, q \leq n$, $1\leq i,j,k,l\leq d$.  Following the notations of Section~\ref{subsec.matrix.Dirichlet.domain},  in this model we have 
\beqnas
(\bfA^{p, q}_{ij, kl})^{ab, cd}_{r, s} &=& \delta_{pq}\delta_{pr}A_{cd}(\delta_{ak}\delta_{bj}\delta_{ci}\delta_{dl} + \delta_{ai}\delta_{bl}\delta_{ck}\delta_{dj}) - \delta_{pr}\delta_{qs}(A_{kj}\delta_{ai}\delta_{bc}\delta_{dl} + A_{il}\delta_{ck}\delta_{ad}\delta_{bj}) \\
& &+ \delta_{pr}\delta_{qs}( B_{ia, ld}\delta_{bj}\delta_{ck} +  B_{bj, ck}\delta_{ai}\delta_{dl} - B_{bj, ld}\delta_{ai}\delta_{ck} - B_{ia, ck}\delta_{bj}\delta_{dl}).
\eeqnas
Then $(\Delta_{n,d}, \Gamma_{\bfA, \bfB}, D_{a_1, \cdots, a_{n+1} })$ is a polynomial model.  Moreover, in this case,
\beqna\label{l.dirichlet2}
  \LL_{\bfA, \bfB, \bfa}(Z^{(p)}_{ij}) &=&  2(a_p - 1 + d)A_{ij} - \sum^{n}_{q =1} (a_q - 1 + d)\big((AZ^{(p)})_{ij} + (Z^{(p)}A)_{ij}\big)\\
\nonumber & & -(a_{n+1} -1) \big((AZ^{(p)})_{ij} + (Z^{(p)}A)_{ij}\big)- 2A_{ij}\tr(Z^{(p)})\\
\nonumber  & &\ \  \ +  \sum_{ab} \big(B_{ia, jb}Z^{(p)}_{ab} +  B_{bj, ai}Z^{(p)}_{ab} - B_{ia, ba}Z^{(p)}_{bj} - B_{bj, ba}Z^{(p)}_{ia} \big),
\eeqna
where $\bfZ^{(n+1)} = \Id - \sum^{n}_{p=1}\bfZ^{(p)}$. 
\ethm

\bpf

First, let us show that  equations~\eqref{gamma.dirichlet2} and~\eqref{l.dirichlet2} define a polynomial model. In fact, for $1 \leq p, q \leq n$,
\beqnas
\Gamma(Z^{(p)}_{ij}, \log \det \bfZ^{(q)})
&=& \sum_{kl}(Z^{(q)})^{-1}_{lk}\Gamma(Z^{(p)}_{ij}, Z^{(q)}_{kl})\\
&=& 2A_{ij}\delta_{pq}  - (AZ^{(p)})_{ij} - (Z^{(p)}A)_{ij},\\
\Gamma(Z^{(p)}_{ij}, \log \det \bfZ^{(n+1)}) 
&=&   - (AZ^{(p)})_{ij} - (Z^{(p)}A)_{ij},
\eeqnas
which is a polynomial model by Proposition \ref{prop.poly.model.gal}.


Direct computations yield
\beqnas
\LL(Z^{(p)}_{ij}) &=& \sum^{n+1}_{q=1} (a_q - 1) \Gamma(Z^{(p)}_{ij}, \log \det Z^{(q)}) + \sum^{n}_{q} \partial_{Z^{(q)}_{kl}}\Gamma(Z^{(p)}_{ij}, Z^{(q)}_{kl}) \\
 &=& 2(a_p - 1 + d)A_{ij}  - \sum^{n}_{q =1} (a_q - 1 + d)\big((AZ^{(p)})_{ij} + (Z^{(p)}A)_{ij}\big) \\
 & &-(a_{n+1} -1) \big((AZ^{(p)})_{ij} + (Z^{(p)}A)_{ij}\big)- 2A_{ij}\tr(Z^{(p)})\\
 & & + \sum_{ab} \big( B_{ia, jb}Z^{(p)}_{ab} +  B_{bj, ai}Z^{(p)}_{ab} - B_{ia, ba}Z^{(p)}_{bj} - B_{bj, ba}Z^{(p)}_{ia} \big).
\eeqnas

Now we prove that if $\bfA$ is a $d \times d$ Hermitian and positive-definite matrix and $\bfB$ is a $d^{2} \times d^{2}$ Hermitian and positive-definite matrix, then $\Gamma_{\bfA, \bfB}$ is elliptic inside the matrix simplex $\Delta_{n, d}$. In fact, consider $\Gamma_{\bfA, \bfB}$ as  a $n \times n$ block matrix, and each block is of size $d^{2} \times d^{2}$, then we may write
\beqnas
\Gamma_{\bfA, \bfB} =  \Gamma_{\bfA} + \Gamma_{\bfB}, 
\eeqnas
where $\Gamma_{\bfA}$ is the block matrix containing $\bfA$ and $\Gamma_{\bfB}$ is the block matrix containing $\bfB$.  


Let $(\Lambda^{1}, \cdots, \Lambda^{n})$ be any sequence of $d \times d$ Hermitian matrices. Then,
\beqnas
& &\sum^{n}_{p, q = 1}\sum_{ijkl} \lambda^{p}_{ij}\bar{\lambda}^{q}_{kl}(\Gamma^{p, q}_{\bfA})_{ij, lk}\\
&=& \sum^{n}_{p, q =1}\sum_{ijkl} \lambda^{p}_{ij}\bar{\lambda}^{q}_{kl}(\delta_{pq}(A_{ik}Z^{(p)}_{lj} + A_{lj}Z^{(p)}_{ik}) - A_{lj}(Z^{(p)}Z^{(q)})_{ik} - A_{ik}(Z^{(q)}Z^{(p)})_{lj})\\
&=& \tr\big(A (\sum^{n}_{p=1}\bar{\Lambda}^{p}\bfZ^{(p)}(\Lambda^{p})^{t}  -  \sum^{n}_{p, q = 1}\bar{\Lambda}^{p}\bfZ^{(p)}\bfZ^{(q)}(\Lambda^{q})^{t}) \big) \\
& &  \ \  +  \tr\big(A (\sum^{n}_{p=1}(\Lambda^{p})^{t}\bfZ^{(p)}\bar{\Lambda}^{p}  -  \sum^{n}_{p, q = 1}(\Lambda^{p})^{t}\bfZ^{(p)}\bfZ^{(q)}\bar{\Lambda}^{q}) \big) .
\eeqnas

Then since $A$ is positive definite, we just need  to prove that 
\beqnas
\sum^{n}_{p=1}\bar{\Lambda}^{p}\bfZ^{(p)}(\Lambda^{p})^{t} - \sum^{n}_{p, q = 1}\bar{\Lambda}^{p}\bfZ^{(p)}\bfZ^{(q)}(\Lambda^{q})^{t}, \ \  \sum^{n}_{p=1}(\Lambda^{p})^{t}\bfZ^{(p)}\bar{\Lambda}^{p}  -  \sum^{n}_{p, q = 1}(\Lambda^{p})^{t}\bfZ^{(p)}\bfZ^{(q)}\bar{\Lambda}^{q},
\eeqnas are non negative definite. For the first one, given a vector $X = (X_1, ..., X_d)$, we have
\beqnas
& &X\big(\sum^{n}_{p=1}\bar{\Lambda}^{p}\bfZ^{(p)}(\Lambda^{p})^{t} - \sum^{n}_{p, q = 1}\bar{\Lambda}^{p}\bfZ^{(p)}\bfZ^{(q)}(\Lambda^{q})^{t}\big)X^{*}\\
&=&  X\bar{\Lambda}(\bfZ- YY^{*})\Lambda^{t}X^{*}\\
&=&  X\bar{\Lambda} \bfZ^{\frac{1}{2}}(\Id_{nd} - (\bfZ^{-\frac{1}{2}}Y)(\bfZ^{-\frac{1}{2}}Y)^{*})(X\bar{\Lambda} \bfZ^{\frac{1}{2}})^{*},
\eeqnas
where $\Lambda$ is a vector of matrices $\Lambda = (\Lambda^{1}, ..., \Lambda^{n})$, $\bfZ = \diag(\bfZ^{(1)}, ..., \bfZ^{(n)})$ and $Y$ is a vector of matrices such that $Y = (\bfZ^{(1)}, ..., \bfZ^{(n)})^{*}$. 
Then by Sylvester determinant theorem, we are able to compute the eigenvalues of $\Id - (\bfZ^{-\frac{1}{2}}Y)(\bfZ^{-\frac{1}{2}}Y)^{*}$,
\beqna\label{sylvester}
& &\det\big(\lambda \Id_{nd} - (\Id_{nd} - (\bfZ^{-\frac{1}{2}}Y)(\bfZ^{-\frac{1}{2}}Y)^{*})\big) \\
\nonumber&=& \det\big((\lambda - 1)\Id_{nd} + (\bfZ^{-\frac{1}{2}}Y)(\bfZ^{-\frac{1}{2}}Y)^{*}\big) \\
\nonumber&=& (\lambda - 1)^{(n-1)d}\det(\lambda \Id - \bfZ^{(n+1)}).
\eeqna
Since $\bfZ^{(n+1)} = \Id - \sum^{n}_{p=1}\bfZ^{(p)}$ is also a non negative-definite Hermitian matrix, the above equation means that the eigenvalues of $\Id - (\bfZ^{-\frac{1}{2}}Y)(\bfZ^{-\frac{1}{2}}Y)^{*}$ are all non negative, indicating that it is a non negative definite matrix, such that
$$
X\big(\sum^{n}_{p=1}\bar{\Lambda}^{p}\bfZ^{(p)}(\Lambda^{p})^{t} - \sum^{n}_{p, q = 1}\bar{\Lambda}^{p}\bfZ^{(p)}\bfZ^{(q)}(\Lambda^{q})^{t}\big)X^{*} \geq 0, 
$$ 
thus $\sum^{n}_{p=1}\bar{\Lambda}^{p}\bfZ^{(p)}(\Lambda^{p})^{t} - \sum^{n}_{p, q = 1}\bar{\Lambda}^{p}\bfZ^{(p)}\bfZ^{(q)}(\Lambda^{q})^{t}$ is non negative definite, so is $\sum^{n}_{p=1}(\Lambda^{p})^{t}\bfZ^{(p)}\bar{\Lambda}^{p}  -  \sum^{n}_{p, q = 1}(\Lambda^{p})^{t}\bfZ^{(p)}\bfZ^{(q)}\bar{\Lambda}^{q}$.
Therefore from the fact that $A$ is a positive-definite Hermitian matrix, we have
\beqnas
\tr\big(A (\sum^{n}_{p = 1}\bar{\Lambda}^{p}\bfZ^{(p)}(\Lambda^{p})^{t}  - \sum^{n}_{p, q = 1}\bar{\Lambda}^{p}\bfZ^{(p)}\bfZ^{(q)}(\Lambda^{q})^{t})\big) &\geq& 0,\\
\tr\big(A (\sum^{n}_{p=1}(\Lambda^{p})^{t}\bfZ^{(p)}\bar{\Lambda}^{p}  -  \sum^{n}_{p, q = 1}(\Lambda^{p})^{t}\bfZ^{(p)}\bfZ^{(q)}\bar{\Lambda}^{q}) \big) &\geq& 0.
\eeqnas
Since 
\beqnas
\tr\big(A (\sum^{n}_{p = 1}\bar{\Lambda}^{p}\bfZ^{(p)}(\Lambda^{p})^{t}  - \sum^{n}_{p, q = 1}\bar{\Lambda}^{p}\bfZ^{(p)}\bfZ^{(q)}(\Lambda^{q})^{t})\big) 
= \tr \big(A\bar{\Lambda}(\bfZ- YY^{*})\Lambda^{t}\big), 
\eeqnas
and from equation~\eqref{sylvester},  we know that  the interior of  $\Delta_{n, d}$, $\bfZ- YY^{*}$ is positive definite. Thus following the proof of Lemma \ref{positive-definite}, we may conclude that if 
$$
\tr\big(A (\sum^{n}_{p = 1}\bar{\Lambda}^{p}\bfZ^{(p)}(\Lambda^{p})^{t}  - \sum^{n}_{p, q = 1}\bar{\Lambda}^{p}\bfZ^{(p)}\bfZ^{(q)}(\Lambda^{q})^{t})\big) = 0,
$$ then $\Lambda = 0$. This implies that $\Gamma_{\bfA}$ is elliptic inside $\Delta_{n, d}$. 

As for $\Gamma_{\bfB}$, notice that
\beqnas
& &\sum^{d}_{i,j,k,l = 1} \lambda^{p}_{ij}\bar{\lambda}^{q}_{kl}(\Gamma^{p, q}_{\bfB})_{ij, \bar{kl}}\\
&=& \sum^{d}_{i,j,k,l = 1} \lambda^{p}_{ij}\bar{\lambda}^{q}_{kl}( B_{ia, kb}Z^{(p)}_{aj}Z^{(q)}_{lb} +  B_{aj, bl}Z^{(p)}_{ia}Z^{(q)}_{bk} - B_{aj, kb}Z^{(p)}_{ia}Z^{(q)}_{lb} - B_{ia, bl}Z^{(p)}_{aj}Z^{(q)}_{bk})\\
&=&  (\Lambda^{p}\bar{\bfZ}^{(p)} - \bar{\bfZ}^{(p)}\Lambda^{p})B(\bar{\bfZ}^{(q)}(\Lambda^{q})^{*}- (\Lambda^{q})^{*}\bar{\bfZ}^{(q)}),
\eeqnas

Let $H^{(p)} = \Lambda^{p}\bar{\bfZ}^{(p)} - \bar{\bfZ}^{(p)}\Lambda^{p}$, since $\bfB$ is positive-definite, Hermitian in the sense that $B_{ij, kl} = \bar{B}_{kl, ij}$, then
\beqnas
\sum^{n}_{p, q = 1}\sum^{d}_{i,j,k,l =1} \lambda^{p}_{ij}\bar{\lambda}^{q}_{kl}(\Gamma^{p, q}_{B})_{ij, lk} &=& (\sum^{n}_{p=1}H^{(p)}) B (\sum^{n}_{q=1}(H^{(q)})^{*}) \geq 0,
\eeqnas
which means $\Gamma_{\bfB}$ is non-negative definite.  Since $\Gamma_{\bfA, \bfB} = \Gamma_{\bfA} + \Gamma_{\bfB}$, we know that $\Gamma_{\bfA, \bfB}$ is elliptic inside $\Delta_{n, d}$. Then we finish the proof.

\epf

In what follows, we show that this model may be constructed as  a projection from complex Wishart processes, which are matrix generalizations of Laguerre processes. We first recall the definition of the complex Wishart distribution,
\bdefi\label{complex.wishart}
A $d \times d$ Hermitian positive definite matrix $W$ is said to have a Wishart distribution with parameters $d$, $r \geq d$, if its distribution has a density   given by
\beq\label{wishart.distribution}
C_{r, d} \det(W)^{r - d}e^{- \frac{1}{2}\tr(W)},
\eeq
where $C_{r, d} =(2^{rd}\pi^{\frac{1}{2}d(d-1)}\Gamma(r)\Gamma(r-1)\cdots \Gamma(r -d +1))^{-1}$ is the normalization constant. 
\edefi

When $r \geq d$ is an integer, this distribution can be derived from the Gaussian distributed complex matrix. Indeed,  if we  consider  a $d \times r$ complex matrix $X$ with its elements being  independent Gaussian centered random variables, then $W = XX^{*}$ has the complex Wishart distribution with parameters $d, r$.

A $d \times d$ complex Wishart process $\{W_t, t \geq 0\}$  is usually defined as a solution to the following stochastic differential equation, 
\beqnas
dW_t = \sqrt{W_t}dB_t +  dB^{*}_t\sqrt{W_t} +  (\alpha W_t + \beta\Id_{d})dt, \ \ \  W_t = W_{0},
\eeqnas
where $\{B_t, t \geq 0\}$ is a $d \times d$ complex Brownian motion, $W_0$ is a $d \times d$ Hermitian matrix.  Wishart processes have been deeply studied, see \cite{Bru, demni3, donati-martin} etc. There exists more general form of Wishart processes that we will not consider here.  In what follows, we extend the construction of Wishart laws from matrix Gaussian ones at the level of processes, exactly as in the scalar case where Laguerre processes may be constructed (with suitable parameters) from Ornstein-Uhlenbeck ones. We will apply the matrix extracting procedure again.

The generator of an Ornstein-Uhlenbeck process on  $N \times N$ complex matrices is given, on the entries  $\{z_{ij}\}$ of a complex matrix $\bfz$,   by
\beqnas
\Gamma(z_{ij}, z_{kl}) &=& 0,\\
\Gamma(z_{ij}, \bar{z}_{kl}) &=& 2\delta_{ik}\delta_{jl},\\
\LL(z_{ij}) &=& -z_{ij}.
\eeqnas

This describes a process on matrices  where the entries are independent complex Ornstein-Uhlenbeck processes. 
Now we start the "matrix-extracting" procedure on $\bfz$, as we did before on $SU(N)$ in Section~\ref{construction2}, i.e., take the first $d$ lines and split the $N$ columns into $(n + 1)$ parts $I_{1}, \cdots, I_{n+1}$, such that $d_{i} = |I_{i}|$, for $1 \leq i \leq (n+1)$ and $N = d_1 + ... + d_{n+1}$, then we have $(n + 1)$ extracted matrices $\bfY^{(1)}$, ..., $\bfY^{(n+1)}$. For $1 \leq p \leq (n+1)$, define $\bfW^{(p)} = \bfY^{(p)}(\bfY^{(p)})^{*}$ with entries $W^{(p)}_{ij}$. 

\bprop
$\{\bfW^{(p)}, 1 \leq p \leq (n+1)\}$ form a family of independent complex Wishart processes,  whose reversible measure respectively given by \eqref{wishart.distribution} with $r_{p}= d_{p}$ for $1 \leq p \leq (n+1)$. Moreover, the image of the complex Gaussian measure through the "matrix-extracting" procedure is a product of complex Wishart distributions.  
\eprop

\bpf
One may check
\beqna
\label{wishart1}\Gamma(W^{(p)}_{ij}, W^{(q)}_{kl}) &=& 2\delta_{pq}(\delta_{jk}W^{(p)}_{il} + \delta_{il}W^{(p)}_{kj}),\\
\label{wishart2}\LL(W^{(p)}_{ij}) &=& 4d_{p}\delta_{ij} - 2W^{(p)}_{ij}.
\eeqna

Let $\rho$ be the density of the reversible measure of $\{\bfW^{(1)}, ..., \bfW^{(n+1)}\}$. Then,
$$
\Gamma(\log \rho, W^{p}_{ij}) = 4(d_{p} - d)\delta_{ij} -2W^{(p)}_{ij},
$$
and we also have 
\beqnas
\Gamma(\log \det(W^{(p)}), W^{(p)}_{ij}) &=& 4\delta_{ij},\\
\Gamma(\tr W^{(p)}, W^{(p)}_{ij}) &=& 4W^{(p)}_{ij},
\eeqnas
therefore,
$$
\rho = C\prod^{n+1}_{p=1}\det(W^{(p)})^{d_p - d}e^{-\frac{1}{2}\tr (\sum^{n+1}_{p=1}W^{(p)})},
$$
which shows that, under the reversible measure,  we have a  family of $d \times d$ independent Wishart matrices $\{\bfW^{(1)}, ..., \bfW^{(n+1)}\}$. 
\epf

We now construct  a process on the complex matrix simplex from independent Wishart processes $(\bfW^{(1)}, \cdots \bfW^{(n+1)})$. As in the scalar case, we obtain a kind of warped product on the set $\bfD \times \Delta_{n,d}$, where $\bfD$ denotes the set of real diagonal matrices with positive diagonal entries.

Let $\bfS = \sum^{n+1}_{p=1}\bfW^{(p)}$.  Since $\bfS$ is a positive-definite Hermitian matrix, we may assume that  it has a spectral decomposition $\bfS = \bfU \bfD^{2} \bfU^{*}$, where $\bfU$ is unitary  and $\bfD = \diag\{\lambda_{1}, \cdots, \lambda_{d}\}$. Observe that $\bfU$ is not uniquely determined, since we may change $\bfU$ into $\bfU \bfP$ where $\bfP = \diag\{e^{i\phi_{1}}, \cdots, e^{i\phi_{d}}\}$, $0 \leq \phi_1, \cdots, \phi_d \leq 2\pi$ , and this amounts to the choice of a phase for the eigenvectors. In this paper, we choose $\bfU$ to be the one that has real elements on its diagonal, such that $\bfU$ is an analytic function of $\bfS$ in the Weyl chamber $\{\lambda_1< \cdots < \lambda_d\}$. This choice will be irrelevant to the construction of the process. Moreover, we introduce $(\bfV^{(1)}, \cdots , \bfV^{(d)})$, whose elements are given by
\beq \label{def.Vp}
V^{(p)}_{ij} = U_{ip}U^{*}_{pj}
\eeq
for $1 \leq i, j, p \leq d$, and we see  that the choice of phase in $\bfU$ has no influence on $\{\bfV^{(p)}\}$. 

Then for $1 \leq i \leq (n+1)$, write
\beq \label{def.Mi.Zi}
\bfM^{(i)} = \bfS^{-\frac{1}{2}}\bfW^{(i)}\bfS^{-\frac{1}{2}}, \ \ \ \bfZ^{(i)} = \bfU^{*}\bfM^{(i)}\bfU, 
\eeq
for which we have the following result.

\bthm\label{th.wishart.dirichlet}
Let  $(\bfW^{1}, \cdots, \bfW^{(n+1)})$ be $(n + 1)$ independent Wishart processes. Then with $\bfD = \diag\{\lambda_{1}, \cdots, \lambda_{d}\}$, $0 \leq \lambda_{1} \leq \cdots \leq \lambda_{d}$   and $\bfZ^{(i)}$ defined as in equation~\eqref{def.Mi.Zi}, $(\bfD, \bfZ^{(1)}, \cdots, \bfZ^{(n)})$ is a Markov  diffusion process,  where $(\bfZ^{(1)}, \cdots, \bfZ^{(n)})$ lives in the matrix simplex $\Delta_{d, n}$. The generator of the process is 
\beqna\label{opt.wishart.dirichlet}
\sum^{d}_{i=1}\big(\partial^{2}_{\lambda_i} +(\frac{2(N - d)+ 1}{\lambda_i} - \lambda_i  + \sum_{j \neq i}\frac{4\lambda_i}{\lambda^{2}_i - \lambda^{2}_j})\partial_{\lambda_{i}}\big) + \LL_{\bfA, \bfB, \bfa},
\eeqna
where 
$\LL_{\bfA, \bfB, \bfa}$ is defined in Theorem~\ref{th.Dirichlet.simplex} with 
\beqnas
A_{ij} &=& 2\lambda_i^{-2}\delta_{ij}, \\  
B_{ij, kl} &=& \left\{ 
\begin{array}{ll}
2\frac{\lambda^{2}_{i} + \lambda^{2}_{j}}{(\lambda^{2}_{i} - \lambda^{2}_{j})^{2}}\delta_{ik}\delta_{jl}, & \hbox{$i \neq j$ and $k \neq l$,}\\
\frac{1}{\lambda^{2}_{i}}\delta_{ij}\delta_{ik}\delta_{jl}, & \hbox{$i = j$ or $k = l$.}
\end{array}
\right.
\eeqnas
for $1 \leq i, j, k, l \leq d$ and $a_p = d_p - d +1$ for $1 \leq p \leq (n+1)$.
\ethm

It is known that starting from random matrices $(W^{(1)}, \cdots, W^{(n+1)})$ distributed as independent Wishart distributions, one could get the matrix Dirichlet distribution through $M^{i} = S^{-\frac{1}{2}}W^{(i)}S^{-\frac{1}{2}}$, where $S = \sum^{n+1}_{i = 1}W^{i}$, see \cite{guptanagar}. Also from our results in the scalar case (Section \ref{sec.scalar.simplex}), it is natural to guess that $(\bfM^{(1)}, \cdots, \bfM^{(n)})$ may be a matrix Dirichlet process. As we will see in the following proposition, given $\bfS$,  $(\bfM^{(1)}, \cdots, \bfM^{(n)})$ is indeed a matrix Dirichlet process; However, the operator of $(\bfS, \bfM^{(1)}, \cdots, \bfM^{(n)})$ is much more complicated than the one of $(\bfD, \bfZ^{(1)}, \cdots, \bfZ^{(n)})$, since $\Gamma(\bfS, \bfM^{(i)}) \neq 0$, and therefore does not have the structure of a (generalized)  warped product.

\bprop\label{prop.m}
$(\bfS, \bfM^{(1)}, \cdots, \bfM^{(n)})$ is a Markov diffusion process where $(\bfM^{(1)}, \cdots, \bfM^{(n)})$ lives on the matrix simplex $\Delta_{n, d}$, and the generator of the process is 
\beqna\label{operator.sm}
 & \sum_{ijkl}2(\delta_{jk}S_{il} + \delta_{il}S_{kj})\partial_{S_{ij}}\partial_{S_{kl}} + \sum_{ij}(4N\delta_{ij} - 2S_{ij})\partial_{S_{ij}}  + \LL_{\bfA, \bfB, \bfa}\\
\nonumber&  \  \ + \sum_{ijkl, p} \sum^{d}_{a, b =1} 2\frac{\lambda_a - \lambda_b}{\lambda_a + \lambda_b}((M^{(p)}V^{(a)})_{il}V^{(b)}_{kj} - V^{(a)}_{il}(V^{(b)}M^{(p)})_{kj})\partial_{S_{ij}}\partial_{M^{(p)}_{kl}},\\
\nonumber&
\eeqna
where $V^{(i)}$ defined by equation~\eqref{def.Vp},  $\LL_{\bfA, \bfB, \bfa}$ is defined in Theorem~\ref{th.Dirichlet.simplex}  with
$$
A = 2S^{-1}, \ \  B_{ij, kl} = \sum^{d}_{r,s=1}\frac{4}{(\lambda_{r} + \lambda_s)^{2}}V^{(r)}_{ik}V^{(s)}_{lj},
$$ and $a_p -1 + d = d_p$.
\eprop

Before proving Theorem \ref{th.wishart.dirichlet} and Proposition \ref{prop.m}, we first give the following lemmas regarding the action of the diffusion operators of $\bfN = \bfS^{\frac{1}{2}}$, $\{\bfM^{(i)}, 1 \leq i \leq n\}$ and $\{\bfZ^{(i)}, 1 \leq i \leq n\}$.

\blem\label{s.eigenvalues}
For the generator of independent Wishart matrices \\
$(\bfW^{(1)}, \cdots, \bfW^{(n+1)})$ and $\bfS = \sum^{n+1}_{p=1}\bfW^{(p)}$, we have
\beqna
\label{zeq1}\Gamma(S_{ij}, S_{kl}) &=&  2(\delta_{jk}S_{il} + \delta_{il}S_{kj}),\\
\label{zeq3}\LL(S_{ij}) &=& 4N\delta_{ij} - 2S_{ij},\\
\label{zeq2}\Gamma(S_{ij}, W^{(p)}_{kl}) &=&   2(\delta_{jk}W^{(p)}_{il} + \delta_{il}W^{(p)}_{kj}).
\eeqna
Moreover, suppose $\bfS$ has a spectral decomposition $\bfS = \bfU \bfD^{2} \bfU^{*}$, where $\bfU$ is unitary and $\bfD = \diag\{\lambda_{1}, \cdots, \lambda_{d}\}$, we have
\beqna
\label{zeq111}\Gamma(\lambda_i, \lambda_j) &=& \delta_{ij}, \\
\label{zeq222}\LL(\lambda_i) &=& \frac{2(N - d)+ 1}{\lambda_i} - \lambda_i  + 4\lambda_i\sum_{j \neq i}\frac{1}{\lambda^{2}_i - \lambda^{2}_j}.
\eeqna
\elem

\bpf
Formulas \eqref{zeq1}, \eqref{zeq2} and \eqref{zeq3} are straight-forward from \eqref{wishart1} and \eqref{wishart2}.

By \eqref{zeq1}, \eqref{zeq3}, we are able to compute the diffusion operators of  $D = \{\lambda_{1}, ..., \lambda_{d}\}$ by the method described  in Theorem~\ref{th.polar} in Appendix~\ref{polar}. Notice that the $\Gamma$ operator of $\bfS$ is the same as $\bfH$ in Theorem \ref{th.polar}, which implies that 
\beqnas
\Gamma(\lambda_i, \lambda_j) = \delta_{ij}. 
\eeqnas
 Now we compute $\LL(\lambda_i)$. First, let $P(X) = \det(S - X\Id)$, then by \eqref{zeq1}, \eqref{zeq3} we have
\beqna\label{zeq16}
& &\LL(\log P(X)) \\
\nonumber&=&  - 4\tr((S - X\Id)^{-1})\tr((S - X\Id)^{-1}S) + 4N\tr((S - X\Id)^{-1})\\
\nonumber& & - 2\tr((S - X\Id)^{-1}S)\\
\nonumber &=&  4(d - N)\frac{P'(X)}{P(X)} - 4\frac{X(P'(X))^{2}}{P^{2}(X)})  - 2d + 2\frac{XP'(X)}{P(X)}).
\eeqna
On the other hand, let $\eta_{i} = \lambda^{2}_{i}$ be the eigenvalues of $\bfS$, we have
\beqnas
\LL(\log P(X)) = \sum^{d}_{i=1}(-\frac{4\eta_i}{(\eta_i - X)^{2}} + \frac{\LL(\eta_{i})}{\eta_i - X}).
\eeqnas
Comparing it with \eqref{zeq16} leads to 
\beqnas
\LL(\eta_{i}) =  8\eta_i\sum_{j \neq i}\frac{1}{\eta_i - \eta_j} + 4(N - d + 1) - 2\eta_i,
\eeqnas
such that 
\beqnas
\LL(\lambda_i) = \frac{2(N - d)+ 1}{\lambda_i} - \lambda_i  + 4\lambda_i\sum_{j \neq i}\frac{1}{\lambda^{2}_i - \lambda^{2}_j}.
\eeqnas
\epf

\blem
Let $S$ be a positive definite, Hermitian matrix and $N$ be its positive definite square root. Suppose $S$ has the spectral decomposition $S = UD^{2}U^{*}$, where $U$ is the unitary part and $D = \diag\{\lambda_1, \cdots, \lambda_d\}$, with $\lambda_i$ all positive.  With $N= UDU^*$,  and  $V^{(p)}_{ij} = U_{ip}U^{*}_{pj}$ for $1 \leq i, j, p \leq d$, we have  for $1 \leq i, j, k, l \leq d$
\beqna\label{sqrt.matrix}
 \partial_{S_{kl}}N_{ij} = \sum^{d}_{r, s = 1}\frac{1}{\lambda_r + \lambda_s}V^{(r)}_{ik}V^{(s)}_{lj}.
\eeqna

\elem

\bpf
For fixed $k, l$, write $\cD^{kl}_{ij} = \partial_{S_{kl}}N_{ij}$. Then from $N^{2}_{ij} = S_{ij}$, we have
\beqnas
\cD^{kl}\bfN + \bfN \cD^{kl} = E^{kl},
\eeqnas
where $E^{kl}$ is the matrix satisfying $E^{kl}_{ab} = \delta_{ak}\delta_{bl}$. Since $\bfS = \bfU \bfD^{2} \bfU^{*}$, we may write $\bfN = \bfU \bfD \bfU^{*}$, such that
\beqnas
(\bfU^{*} \cD^{kl} \bfU) \bfD + \bfD(\bfU^{*} \cD^{kl} \bfU) = \bfU^{*}E^{kl}\bfU.
\eeqnas
Therefore,
\beqnas
(U^{*}\cD^{kl}U)_{ij} = \frac{U^{*}_{ik}U_{lj}}{\lambda_i + \lambda_j},
\eeqnas
so that
\beqnas
\cD^{kl}_{ij} &=& \sum_{r, s}\frac{1}{\lambda_r + \lambda_s}U_{ir}U^{*}_{sj}U^{*}_{rk}U_{ls}\\
&=& \sum_{r, s}\frac{1}{\lambda_r + \lambda_s}V^{(r)}_{ik}V^{(s)}_{lj},
\eeqnas
which finishes the proof. 
\epf

\blem
The following formulas hold for $\bfN = \bfS^{\frac{1}{2}}$, 
\beqna
 \label{zeq4}\Gamma(N_{ij}, W^{(p)}_{kl}) &=& \sum^{d}_{r, s = 1}\frac{2}{\lambda_r + \lambda_s}V^{(r)}_{il}(W^{(p)}V^{(s)})_{kj} +  \sum^{d}_{r, s =1}\frac{2}{\lambda_r + \lambda_s}V^{(s)}_{kj}(V^{(r)}W^{(p)})_{il},\\
\label{zeq5}\Gamma(N_{ij}, N_{kl}) &=& \sum^{d}_{r, s = 1} 2\frac{\lambda^{2}_{r} + \lambda^{2}_{s}}{(\lambda_r + \lambda_s)^{2}}V^{(r)}_{il}V^{(s)}_{kj},\\
\label{zeq6}\LL(N_{ij}) &=& 4\sum^{d}_{r, s=1} \frac{\lambda_s}{(\lambda_r + \lambda_s)^{2}}V^{(r)}_{ij} - N_{ij} + 2(N - d)N^{-1}_{ij}.
\eeqna
Furthermore,
\beqna
\label{zeq7} \Gamma(N^{-1}_{ij}, N_{kl}) &=& -\sum^{d}_{r, s = 1} 2\frac{\lambda^{2}_{r} + \lambda^{2}_{s}}{\lambda_r\lambda_s(\lambda_r + \lambda_s)^{2}}V^{(r)}_{il}V^{(s)}_{kj},\\
\label{zeq8} \Gamma(N^{-1}_{ij}, W^{(p)}_{kl}) &=& - 2\sum^{d}_{r, s =1}\frac{1}{\lambda_r\lambda_s(\lambda_r + \lambda_s)}(V^{(r)}_{il}(W^{(p)}V^{(s)})_{kj} + V^{(s)}_{kj}(V^{(r)}W^{(p)})_{il}),\\
\label{zeq9} \Gamma(N^{-1}_{ij}, N^{-1}_{kl}) &=& 2\sum^{d}_{r, s =1} \frac{\lambda^{2}_{r} + \lambda^{2}_{s}}{\lambda^{2}_r\lambda^{2}_s(\lambda_r + \lambda_s)^{2}}V^{(r)}_{il}V^{(s)}_{kj},\\
\label{zeq10} \LL(N^{-1}_{ij})  &=& 4\sum^{d}_{r, s =1} \frac{1}{\lambda_s(\lambda_r + \lambda_s)^{2}}V^{(r)}_{ij} + N^{-1}_{ij} - 2(N- d)(S^{-1}N^{-1})_{ij}.
\eeqna
\elem

\bpf

By \eqref{sqrt.matrix},  we are able to compute
\beqnas
\Gamma(N_{ij}, W^{(p)}_{kl}) &=& \sum_{ab}\cD^{ab}_{ij}\Gamma(S_{ab}, W^{(p)}_{kl})\\
&=& 2\sum^{d}_{r, s =1}\frac{1}{\lambda_r + \lambda_s}V^{(r)}_{il}(W^{(p)}V^{(s)})_{kj} +  2\sum^{d}_{r, s =1}\frac{1}{\lambda_r + \lambda_s}V^{(s)}_{kj}(V^{(r)}W^{(p)})_{il},
\eeqnas

and 
\beqnas
\Gamma(N_{ij}, N_{kl}) &=& \sum_{abcd} \cD^{ab}_{ij}\cD^{cd}_{kl}\Gamma(S_{ab}, S_{cd})\\
&=& 2\sum^{d}_{r, s=1} \frac{\lambda^{2}_{r} + \lambda^{2}_{s}}{(\lambda_r + \lambda_s)^{2}}V^{(r)}_{il}V^{(s)}_{kj},\\
\LL(N_{ij}) 
&=&  4\sum \frac{\lambda_s}{(\lambda_r + \lambda_s)^{2}}V^{(r)}_{ij} - N_{ij} + 2(N - d)N^{-1}_{ij}.
\eeqnas


Moreover, due to the fact that
$$
\partial_{N_{ab}}N^{-1}_{ij} = - N^{-1}_{ia}N^{-1}_{bj}, 
$$
we have
\beqnas
& &\Gamma(N^{-1}_{ij}, W^{(p)}_{kl}) = \sum^{d}_{a, b=1}\partial_{N_{ab}}N^{-1}_{ij}\Gamma(N_{ab}, W^{(p)}_{kl})\\
  &=&  - 2\sum^{d}_{r, s =1}\frac{1}{\lambda_r\lambda_s(\lambda_r + \lambda_s)}V^{(r)}_{il}(W^{(p)}V^{(s)})_{kj} - 2\sum^{d}_{r, s=1}\frac{1}{\lambda_r\lambda_s(\lambda_r + \lambda_s)}V^{(s)}_{kj}(V^{(r)}W^{(p)})_{il}.
 \eeqnas
 
 Similar computations yield 
 \beqnas
\Gamma(N^{-1}_{ij}, N_{kl}) 
&=& -2\sum^{d}_{r, s=1} \frac{\lambda^{2}_{r} + \lambda^{2}_{s}}{\lambda_r\lambda_s(\lambda_r + \lambda_s)^{2}}V^{(r)}_{il}V^{(s)}_{kj},\\
\Gamma(N^{-1}_{ij}, N^{-1}_{kl})
 &=& 2\sum^{d}_{r, s=1} \frac{\lambda^{2}_{r} + \lambda^{2}_{s}}{\lambda^{2}_r\lambda^{2}_s(\lambda_r + \lambda_s)^{2}}V^{(r)}_{il}V^{(s)}_{kj},\\
 \LL(N^{-1}_{ij}) 
 &=& \sum^{d}_{k, l=1}\partial_{N_{kl}}N^{-1}_{ij}\LL(N_{kl}) + \sum^{d}_{k, l, a, b=1}\Gamma(N_{kl}, N_{ab})\partial_{N_{ab}}\partial_{N_{kl}}N^{-1}_{ij}\\
 &=& 4\sum^{d}_{r,s=1} \frac{1}{\lambda_s(\lambda_r + \lambda_s)^{2}}V^{(r)}_{ij} + N^{-1}_{ij} - 2(N - d)(S^{-1}N^{-1})_{ij}.
\eeqnas 

\epf

\blem\label{lemma.m}
With $V^{(i)}$ defined in equation~\eqref{def.Vp},  we have 
\beqna\label{zeq11}
& &\Gamma(M^{(p)}_{ij}, M^{(q)}_{kl})\\
\nonumber&=& 2\delta_{pq}(S^{-1}_{il}M^{(p)}_{kj} + S^{-1}_{kj}M^{(p)}_{il}) - 2S^{-1}_{kj}(M^{(p)}M^{(q)})_{il} - 2S^{-1}_{il}(M^{(q)}M^{(p)})_{kj}\\
\nonumber& & - \sum^{d}_{a, b=1} \frac{4}{(\lambda_a + \lambda_b)^{2}}(V^{(b)}M^{(p)})_{kj}(V^{(a)}M^{(q)})_{il}  - \sum^{d}_{a, b=1} \frac{4}{(\lambda_a + \lambda_b)^{2}}(M^{(p)}V^{(a)})_{il}(M^{(q)}V^{(b)})_{kj} \\
\nonumber& & +\sum^{d}_{a, b=1} \frac{4}{(\lambda_a + \lambda_b)^{2}}V^{(a)}_{kj}(M^{(p)}V^{(b)}M^{(q)})_{il}   +\sum^{d}_{a, b =1}\frac{4}{(\lambda_a + \lambda_b)^{2}}V^{(a)}_{il}(M^{(q)}V^{(b)}M^{(p)})_{kj},
\eeqna
and
\beqna\label{zeq12}
& &\LL(M^{(p)}_{ij}) \\
\nonumber &=&  4d_{p}S^{-1}_{ij} - 2(N - d)(S^{-1}M^{(p)})_{ij} - 2(N -d)(M^{(p)}S^{-1})_{ij}- 4S^{-1}_{ij}\tr(M^{(p)})\\
\nonumber & &- 4\sum^{d}_{a, b=1}\frac{1}{(\lambda_a + \lambda_b)^{2}}(V^{(a)}M^{(p)})_{ij} - 4\sum^{d}_{a, b=1}\frac{1}{(\lambda_a + \lambda_b)^{2}}(M^{(p)}V^{(b)})_{ij}\\
\nonumber & & + 8\sum^{d}_{a, b=1}\frac{1}{(\lambda_a + \lambda_b)^{2}}V^{(a)}_{ij}\tr(V^{(b)}M^{(p)}).
\eeqna
Moreover, we have
\beqna
 \label{zeq20}\Gamma(M^{(p)}_{ij}, S_{kl}) &=& \sum^{d}_{a, b =1} 2\frac{\lambda_a - \lambda_b}{\lambda_a + \lambda_b}((M^{(p)}V^{(a)})_{il}V^{(b)}_{kj} - V^{(a)}_{il}(V^{(b)}M^{(p)})_{kj}).
\eeqna
\elem

\bpf
Since 
\beqnas
\Gamma(M^{(p)}_{ij}, M^{(q)}_{kl}) =  \Gamma((N^{-1}W^{(p)}N^{-1})_{ij}, (N^{-1}W^{(q)}N^{-1})_{kl}),
\eeqnas
then by \eqref{wishart1}, \eqref{zeq8} and \eqref{zeq9} we are able to prove $\eqref{zeq11}$.

As for \eqref{zeq12}, direct computations yield
\beqnas
& &\LL(M^{(p)}_{ij}) = \sum^{d}_{r,s=1}\LL(N^{-1}_{ir}W^{(p)}_{rs}N^{-1}_{sj})\\
&=& \sum^{d}_{r,s=1} \big( 2\Gamma(N^{-1}_{ir}, N^{-1}_{sj})W^{p}_{rs} + 2\Gamma(N^{-1}_{ir}, W^{p}_{rs})N^{-1}_{sj} + 2N^{-1}_{ir}\Gamma(W^{p}_{rs}, N^{-1}_{sj})\\
& & + \LL(N^{-1}_{ir})W^{p}_{rs}N^{-1}_{sj} + N^{-1}_{ir}\LL(W^{p}_{rs})N^{-1}_{sj} + N^{-1}_{ir}W^{p}_{rs}\LL(N^{-1}_{sj}) \big)\\
&=&  8\sum^{d}_{a, b=1}\frac{1}{(\lambda_a + \lambda_b)^{2}}V^{(a)}_{ij}\tr(V^{(b)}M^{(p)}) - 4S^{-1}_{ij}\tr(M^{p})\\
& &- 4\sum^{d}_{a, b=1}\frac{1}{(\lambda_a + \lambda_b)^{2}}(V^{(a)}M^{(p)})_{ij} - 4\sum^{d}_{a, b=1}\frac{1}{(\lambda_a + \lambda_b)^{2}}(M^{(p)}V^{(b)})_{ij} \\
& &  + 4d_{p}S^{-1}_{ij} - 2(N - d)(S^{-1}M^{(p)})_{ij} - 2(N-d)(M^{(p)}S^{-1})_{ij},
\eeqnas
where the last equality is due to \eqref{wishart2}, \eqref{zeq8}, \eqref{zeq9} and \eqref{zeq10}.

By \eqref{zeq4}, \eqref{zeq7}, we get
\beqna\label{zeq14444}
& &\Gamma(M^{(p)}_{ij}, N_{kl}) \\
\nonumber&=& \sum^{d}_{a =1}\Gamma(N_{kl}, N^{-1}_{ia})(W^{(p)}N^{-1})_{aj} + \sum^{d}_{a, b =1}N^{-1}_{ia}\Gamma(N_{kl}, W^{(p)}_{ab})N^{-1}_{bj} \\
\nonumber& & \ \ \ \ + \sum^{d}_{a =1}(N^{-1}W^{(p)})_{ia}\Gamma(N_{kl}, N^{-1}_{aj})\\
\nonumber&=& \sum^{d}_{a, b=1}2\frac{\lambda_a - \lambda_b}{(\lambda_a + \lambda_b)^{2}}(M^{(p)}V^{(a)})_{il}V^{(b)}_{kj} -  2\sum^{d}_{a, b=1}\frac{\lambda_a - \lambda_b}{(\lambda_a + \lambda_b)^{2}}V^{(a)}_{il}(V^{(b)}M^{(p)})_{kj},
\eeqna

which leads to \eqref{zeq20},
showing that in fact $\bfM^{(p)}$ and $\bfS$ are not independent.  
\epf

Now we may obtain Proposition \ref{prop.m} directly from Lemma \ref{lemma.m}.

\bpf
Combining  \eqref{zeq1}, \eqref{zeq3} and \eqref{zeq11}, \eqref{zeq12} and \eqref{zeq20}, we derive  \eqref{operator.sm} as the operator of $(\bfS, \bfM^{(1)}, \cdots, \bfM^{(n)})$.
\epf

\brmq
It is natural to consider the reversible measure of  \\ $(\bfS, \bfM^{(1)}, \cdots, \bfM^{(n)})$ as soon as we have Proposition \ref{prop.m}. However, it appears quite complicated in this case, due to the fact that we need to take the unitary part $\bfU$ into account. More precisely, we define $\rho$ to be the density of the reversible measure of $(\bfS, \bfM^{(1)}, \cdots, \bfM^{(n)})$ and decompose $\bfS$ into $(\bfD, \bfU)$, then we may write $\rho = \rho_{\bfD}\rho_{\bfU}D_{\bbC, \bfa}$, where $\rho_{\bfD}$ is obtained from the generators of $\{\lambda_{i}\}$, see Theorem \ref{th.wishart.dirichlet}. Again by \eqref{eq.density}, we obtain the formulas of $\rho_{\bfU}$, which are so complicated that we are not able to provide any explicit expression.

\ermq

To examine the relation of $\bfN$ and $\bfM^{(p)}$ more precisely, we decompose $\bfN$ into $\bfD$ and $\bfU$, and explore their relations with $\bfM^{(p)}$ separately. The method is adapted from \cite{BakryZ}, i.e., to obtain the action of the  $\Gamma$ operator  of the spectrum of $\bfN$ and $\bfM^{(p)}$ by computing its action  on the characteristic polynomial of $\bfN$ and $\bfM^{(p)}$.

\blem\label{lemma.wishart.eq.Gamma.L}
\beqna\label{zeq1444}
\Gamma(M^{(p)}_{ij}, \lambda_k) &=& 0,
\eeqna
and
\beqna
\label{zeq14}\Gamma(M^{(p)}_{ij}, U_{kl}) &=&  \sum^{d}_{a=1}g_{al}(M^{(p)}U)_{il}U^{*}_{aj}U_{ka}  -  \sum^{d}_{a=1}g_{al}U_{il}(U^{*}M^{(p)})_{aj}U_{ka},\\
\label{zeq144}\Gamma(U^{*}_{kl}, M^{(p)}_{ij}) &=& - \sum^{d}_{a=1}g_{ka}(M^{(p)}U)_{ia}U^{*}_{al}U^{*}_{kj} + \sum^{d}_{a=1}g_{ka}U_{ia}U^{*}_{al}(U^{*}M^{(p)})_{kj},
\eeqna
where
\beqnas
g_{ij} = \left\{ 
\begin{array}{ll}
\frac{2}{(\lambda_i + \lambda_j)^{2}}, & \hbox{$i \neq j$,}\\
0, & \hbox{$i = j$,}
\end{array}
\right.
\eeqnas
for $1 \leq i, j \leq d$.
\elem

\bpf
Let $P_{\bfN}(X) = \det(\bfN -  X\Id) = \prod^{d}_{i=1}(\lambda_i - X)$ be the characteristic polynomial of $\bfN$. Notice that 
\beqnas
& &\Gamma(\log P_{N}(X), M^{(p)}_{ij}) = \sum^{d}_{k, l=1} N^{-1}(X)_{lk}\Gamma(N_{kl}, M^{(p)}_{ij})\\
&=&  \sum_{r, s}2\frac{\lambda_r - \lambda_s}{(\lambda_r + \lambda_s)^{2}}\frac{1}{\lambda_r - X}\delta_{rs}U^{*}_{sj}(M^{(p)}U)_{ir} -  2\sum_{r, s}\frac{\lambda_r - \lambda_s}{(\lambda_r + \lambda_s)^{2}}\frac{1}{\lambda_r - X}\delta_{rs}U_{ir}(U^{*}M^{(p)})_{sj} \\
&=& 0,
\eeqnas
from which we deduce that 
$$
\Gamma(M^{(p)}_{ij}, \lambda_k) = 0,
$$
since 
$$
\Gamma(\log P_{N}(X), M^{(p)}_{ij}) = \sum^{d}_{k =1} \frac{1}{\lambda_{k} - X}\Gamma(\lambda_k, M^{(p)}_{ij}) = 0.
$$

Also notice that $\Gamma(N_{kl}, M^{(p)}_{ij})$ is invariant under the unitary transformation $$(\bfN, \bfM^{(p)}) \rightarrow (\bfU^{0}\bfN(\bfU^{0})^{*}, \bfU^{0}\bfM^{(p)}(\bfU^{0})^{*})$$ for any unitary matrix $U^{0}$, (for the details of the invariance of transformation of diffusion operators, see Appendix \ref{polar}), we may compute $\Gamma(M^{(p)}_{ij}, U_{kl})$ at $U = \Id$, then obtain it at any $U$. More precisely, by \eqref{zeq14444}, we obtain at $U = \Id$, 
$$
\Gamma(N_{kl}, M^{(p)}_{ij}) = 2\frac{\lambda_l - \lambda_k}{(\lambda_l + \lambda_k)^{2}}\delta_{kj}M^{(p)}_{il} -  2\frac{\lambda_l - \lambda_k}{(\lambda_l + \lambda_k)^{2}}\delta_{il}M^{(p)}_{kj},
$$

and on the other hand at $U=\Id$, 
$$
\Gamma(N_{kl}, M^{(p)}_{ij}) = \Gamma(U_{kl}, M^{(p)}_{ij})(\lambda_{l} - \lambda_{k}).
$$

Combining the two equalities together, we have for $k \neq l$
\beqna\label{zeq15}
\Gamma(U_{kl}, M^{(p)}_{ij})(U = \Id) =  2\frac{1}{(\lambda_l + \lambda_k)^{2}}(\delta_{kj}M^{(p)}_{il} -  \delta_{il}M^{(p)}_{kj}).
\eeqna

Moreover, at $U = \Id$,  for $1 \leq k \leq d$ 
$$
\Gamma(U_{kk}, \cdot) =  - \Gamma(U^{*}_{kk}, \cdot). 
$$
Then by the fact that $U_{kk}$ is real,  we have at $U = \Id$,
$$
\Gamma(U_{kk}, M^{(p)}_{ij}) = 0.
$$

Now write 
\beqnas
g_{ij} = \left\{ 
\begin{array}{ll}
\frac{2}{(\lambda_i + \lambda_j)^{2}}, & \hbox{$i \neq j$,}\\
0, & \hbox{$i = j$,}
\end{array}
\right.
\eeqnas

then at any $U$, we derive \eqref{zeq14} and \eqref{zeq144} by the invariance under unitary transformations.



\epf



\blem\label{lm.dirichlet2m}
The diffusion operators of $\{\bfZ^{i}\}$ are given by 
\beqna\label{gamma.dirichlet2m}
& &\Gamma(Z^{(p)}_{ij}, Z^{(q)}_{kl})\\
\nonumber&=& 2\delta_{pq}(\delta_{il}\frac{1}{\lambda^{2}_{i}}Z^{(p)}_{kj} + \delta_{kj}\frac{1}{\lambda^{2}_{j}}Z^{(p)}_{il}) - 2(\delta_{kj}\frac{1}{\lambda^{2}_{j}}(Z^{(p)}Z^{(q)})_{il} + \delta_{il}\frac{1}{\lambda^{2}_{i}}(Z^{(q)}Z^{(p)})_{kj})\\
\nonumber& & +  \sum^{d}_{a=1}(y_{ia}\delta_{il}Z^{(p)}_{aj}Z^{(q)}_{ka} + y_{ka}\delta_{kj}Z^{(p)}_{ia}Z^{(q)}_{al}) - y_{ik}Z^{(p)}_{kj}Z^{(q)}_{il} - y_{jl}Z^{(p)}_{il}Z^{(q)}_{kj},
\eeqna
where $y_{ij} = 2\frac{\lambda^{2}_{i} + \lambda^{2}_{j}}{(\lambda^{2}_{i} - \lambda^{2}_{j})^{2}}$ when $i \neq j$ and $y_{ii} = \frac{1}{\lambda_{i}^{2}}$ for $1 \leq i \leq d$,
and 
\beqna\label{l.dirichlet2m}
 \LL(Z^{(p)}_{ij}) &=& 4d_{p}\frac{1}{\lambda^{2}_{i}}\delta_{ij} - 2(N - d)(\frac{1}{\lambda^{2}_{i}} + \frac{1}{\lambda^{2}_{j}})Z^{(p)}_{ij} - 4\frac{1}{\lambda^{2}_{i}}\delta_{ij}\tr(Z^{(p)})\\
\nonumber& & + 2\sum^{d}_{a=1}y_{ia}\delta_{ij}Z^{(p)}_{aa}  - \sum^{d}_{a=1}y_{ja}Z^{(p)}_{ij} - \sum^{d}_{a =1}y_{ia}Z^{(p)}_{ij}.
\eeqna
Moreover, 
\beqna
\label{z.diag}\Gamma(Z^{(p)}_{ij}, \lambda_k) &=& 0,\\
\label{z.unitary}\Gamma(Z^{(p)}_{ij}, U_{kl})  &=&\delta_{il}\sum^{d}_{a \neq l}d_{al}U_{ka}Z^{(p)}_{aj} - d_{jl}U_{kj}Z^{(p)}_{il},
\eeqna
where $d_{ij} = 4\frac{\lambda_i\lambda_j}{(\lambda^{2}_i - \lambda^{2}_j)^{2}}$ when $i \neq j$ and $d_{ii} = 0$.
\elem

\bpf
Notice the diffusion generator of $(\bfU, \bfM^{(1)}, \cdots, \bfM^{(n)})$ is  invariant through the map
$$
(\bfU, \bfM^{(1)}, \cdots, \bfM^{(n)}) \rightarrow (\bfU^{0}\bfU, \bfU^{0}\bfM^{(1)}(\bfU^{0})^{*}, \cdots, \bfU^{0}\bfM^{(n)}(\bfU^{0})^{*}).
$$ Hence the diffusion generator  of $\bfZ = (\bfZ^{(1)}, \cdots, \bfZ^{(n)})$ is also invariant. Therefore to compute $\Gamma(\bfZ^{(p)}, \bfZ^{(q)})$, we may first consider the case at $\bfU = \Id$.
By direct computations, we have
\beqnas
\Gamma(Z^{(p)}_{ij}, Z^{(q)}_{kl}) 
&=& \sum^{d}_{r, s =1} \big(\Gamma(M^{(p)}_{ij}, M^{(q)}_{kl}) + \Gamma(M^{(p)}_{ij}, U^{*}_{kr}U_{sl})M^{(q)}_{rs} + M^{(p)}_{rs}\Gamma(U^{*}_{ir}U_{sj}, M^{(q)}_{kl})\big)\\
& & + \sum^{d}_{u,v,r,s =1}\Gamma(U^{*}_{iu}U_{vj}, U^{*}_{kr}U_{sl})M^{(p)}_{uv}M^{(q)}_{rs} .
\eeqnas

The first term $\Gamma(M^{(p)}_{ij}, M^{(q)}_{kl})$ at $U = \Id$ is straightforward from \eqref{zeq11}.

By \eqref{zeq15}, we get
\beqnas
& &\sum^{d}_{r,s=1}\Gamma(M^{(p)}_{ij}, U^{*}_{kr}U_{sl})M^{(q)}_{rs} \\
&=& \sum^{d}_{a=1}(- g_{ka}\delta_{kj}Z^{(p)}_{ia}Z^{(q)}_{al} -  g_{al}\delta_{il}Z^{(q)}_{ka}Z^{(p)}_{aj})  + g_{jl}Z^{(p)}_{il}Z^{(q)}_{kj}  + g_{ik}Z^{(p)}_{kj}Z^{(q)}_{il},
\eeqnas

\beqnas
& &\sum^{d}_{r,s=1}\Gamma(M^{(q)}_{kl}, U^{*}_{ir}U_{sj})M^{(p)}_{rs} \\
&=& \sum^{d}_{a=1}(- g_{ia}\delta_{il}Z^{(q)}_{ka}Z^{(p)}_{aj} -  g_{aj}\delta_{kj}Z^{(p)}_{ia}Z^{(q)}_{al})  + g_{jl}Z^{(p)}_{il}Z^{(q)}_{kj}  + g_{ik}Z^{(p)}_{kj}Z^{(q)}_{il},
\eeqnas
and
\beqnas
& &\sum^{d}_{u,v,r,s =1}\Gamma(U^{*}_{iu}U_{vj}, U^{*}_{kr}U_{sl})M^{(p)}_{uv}M^{(q)}_{rs}\\
&=& \sum^{d}_{a=1}(r_{ia}\delta_{il}M^{(p)}_{aj}M^{(q)}_{ka} + r_{ka}\delta_{kj}M^{(p)}_{ia}M^{(q)}_{al}) - r_{ik}M^{(p)}_{kj}M^{(q)}_{il} - r_{jl}M^{(p)}_{il}M^{(q)}_{kj}.
\eeqnas

Therefore,  collecting all the four terms together we obtain at $U=\Id$,
\beqnas
& &\Gamma(Z^{(p)}_{ij}, Z^{(q)}_{kl})\\
&=&  2\delta_{pq}(\delta_{il}\frac{1}{\lambda^{2}_{i}}Z^{(p)}_{kj} + \delta_{kj}\frac{1}{\lambda^{2}_{j}}Z^{(p)}_{il}) - 2(\delta_{kj}\frac{1}{\lambda^{2}_{j}}(Z^{(p)}Z^{(q)})_{il} + \delta_{il}\frac{1}{\lambda^{2}_{i}}(Z^{(q)}Z^{(p)})_{kj})\\
& & -  \frac{1}{\lambda_{i}^{2}}\delta_{ik}Z^{(p)}_{ij}Z^{(q)}_{il}  - \frac{1}{\lambda_{j}^{2}}\delta_{jl}Z^{(p)}_{ij}Z^{(q)}_{kj}  + \frac{1}{\lambda_{j}^{2}}\delta_{kj}Z^{(p)}_{ij}Z^{(q)}_{jl}   + \frac{1}{\lambda_{i}^{2}}\delta_{il}Z^{(p)}_{ij}Z^{(q)}_{ki}\\
& & + \sum^{d}_{a=1}(r_{ia}\delta_{il}Z^{(p)}_{aj}Z^{(q)}_{ka} + r_{ka}\delta_{kj}Z^{(p)}_{ia}Z^{(q)}_{al}) - r_{ik}Z^{(p)}_{kj}Z^{(q)}_{il} - r_{jl}Z^{(p)}_{il}Z^{(q)}_{kj}.
\eeqnas
 Defining
\beqnas
y_{ij} = \left\{ 
\begin{array}{ll}
r_{ij}, & \hbox{$i \neq j$,}\\
\frac{1}{\lambda_{i}^{2}}, & \hbox{$i = j$,}
\end{array}
\right.
\eeqnas
for $1 \leq i, j \leq d$, we obtain
\beqnas
& &\Gamma(Z^{(p)}_{ij}, Z^{(q)}_{kl})\\
&=& 2\delta_{pq}(\delta_{il}\frac{1}{\lambda^{2}_{i}}Z^{(p)}_{kj} + \delta_{kj}\frac{1}{\lambda^{2}_{j}}Z^{(p)}_{il}) - 2(\delta_{kj}\frac{1}{\lambda^{2}_{j}}(Z^{(p)}Z^{(q)})_{il} + \delta_{il}\frac{1}{\lambda^{2}_{i}}(Z^{(q)}Z^{(p)})_{kj})\\
& & + \sum^{d}_{a=1}(y_{ia}\delta_{il}Z^{(p)}_{aj}Z^{(q)}_{ka} + y_{ka}\delta_{kj}Z^{(p)}_{ia}Z^{(q)}_{al}) - y_{ik}Z^{(p)}_{kj}Z^{(q)}_{il} - y_{jl}Z^{(p)}_{il}Z^{(q)}_{kj}.
\eeqnas
This formula is also valid at any $U$, because $Z^{(p)}$ is invariant under the map $(U, M^{p}) \rightarrow (U^{0}U, U^{0}M^{p}(U^{0})^{*})$.  Thus we derive \eqref{gamma.dirichlet2m}.


By \eqref{zeq15} and Theorem \ref{th.polar} (the diffusion operators of the unitary part $\bfU$ here are the same as in the polar decomposition case), we obtain formula~\eqref{l.dirichlet2m}.

According to Theorem \ref{th.polar}, we know that $\Gamma(U, D) = 0$, together with \eqref{zeq15}, indicating that 
$$
\Gamma(Z^{(p)}_{ij}, \lambda_k) = 0.
$$

Also by Theorem \ref{th.polar},  we have
at $U=\Id$ and at any $U$,
\beqnas
\Gamma(Z^{(p)}_{ij}, U_{kl}) 
=  4\frac{\lambda_k\lambda_l}{(\lambda^{2}_k - \lambda^{2}_l)^{2}}(\delta_{il}Z^{(p)}_{kj} - \delta_{kj}Z^{(p)}_{il}).
\eeqnas
Let $d_{ij} = 4\frac{\lambda_i\lambda_j}{(\lambda^{2}_i - \lambda^{2}_j)^{2}}$ when $i \neq j$ and $d_{ii} = 0$, then at any $U$, we obtain \eqref{z.unitary}.

\epf

Finally we are in the position to prove Theorem \ref{th.wishart.dirichlet}.

\bpf
By \eqref{zeq111}, \eqref{zeq222}, we obtain the generator of $\bfD$. 

Together with \eqref{z.diag}, and comparing \eqref{gamma.dirichlet2m}, \eqref{l.dirichlet2m} with Theorem \ref{model2}, we obtain the operator \eqref{opt.wishart.dirichlet} with
\beqnas
A = 2D^{-2}, \ \ \ B_{ij, kl} = y_{ij}\delta_{ik}\delta_{jl}, 
\eeqnas
and $a_{p} = d_{p} - d + 1$.
\epf

\appendix
\section{Appendix: Polar decomposition of complex matrices}\label{polar}
In this Appendix we aim at describing the polar decomposition of a Brownian 
 complex matrix. The generator of  a Brownian motion on  $d \times d$ complex matrices $\bfm$, with $\{m_{ij}\}$ as its entries, is given by
$$
\Gamma(m_{ij}, m_{kl}) = 0, \ \ \Gamma(m_{ij}, \bar{m}_{kl}) = 2\delta_{ik}\delta_{jl}, \ \ \ \LL(m_{ij}) = 0.
$$

Suppose $\bfm$ has the polar decomposition $\bfm = \bfV \bfN$, where $\bfV$ is unitary and $\bfN =\sqrt{\bfm^{*}\bfm}$ is Hermitian, positive definite.

Suppose $\bfH = \bfm^{*}\bfm$ has the spectral decomposition $\bfH = \bfU\bfD^{2}\bfU^{*}$, where $\bfU$ is unitary and $\bfD = \diag\{x_1, \cdots, x_d\}$ is the diagonal part. Then $\bfN = \bfU \bfD \bfU^{*}$.

Our goal is to describe the generators of  $\bfV$, $\bfN$ and the spectrum part $\bfD$. As for the spectrum part, we use  the method introduced in \cite{BakryZ}, i.e., we consider the generator of the characteristic polynomial. Moreover, our results are based on an important observation that the generator of $\bfm$ is invariant under both left and right unitary multiplications, which leads to the facts that the generators of $\bfH$,  $\bfN$ are invariant under both left and right unitary multiplications, and also that  the generators of $\bfV$ and $\bfU$ are invariant under left unitary multiplications. More precisely,  in our case, the following lemma holds, ensuring that if we know the diffusion operators of $\bfm$, $\bfH$ or $\bfN$  at $\bfV = \bfU = \Id$, then we know them at any $(\bfV, \bfU)$.

\blem\label{invariance}
Let $\bfm, \bfH, \bfN, \bfV, \bfU$ be defined as above and consider a Brownian motion on $\bfm$. Then given unitary matrices $(V, U)$, the generators of $\bfV$ and $\bfU$ are invariant under the map $(\bfV, \bfU) \rightarrow (VU\bfV U^{*}, U\bfU)$. More precisely, the following formulas hold,
\beqna
\label{invariance1}\Gamma(N_{ij}, N_{kl})(V, U) &=& \sum_{pqrs} U_{ip}U^{*}_{qj}U_{kr}U^{*}_{sl}\Gamma(N_{pq}, N_{rs})(\Id, \Id),\\
\label{invariance2}\LL(N_{ij})(V, U) &=& \sum_{pq}U_{ip}U^{*}_{qj}\LL(N_{pq})(\Id, \Id),\\
\label{invariance3}\Gamma(U_{ij}, U_{kl})(V, U) &=& \sum_{pr} U_{ip}U_{kr}\Gamma(U_{pj}, U_{rl})(\Id, \Id),\\
\label{invariance4}\LL(U_{ij})(V, U) &=& \sum_{p} U_{ip}\LL(U_{pj})(\Id, \Id),\\
\label{invariance5}\Gamma(V_{ij}, V_{kl})(V, U) &=& \sum_{pqrs} (VU)_{ip}U^{*}_{qj}(VU)_{kr}U^{*}_{sl}\Gamma(V_{pq}, V_{rs})(\Id, \Id),\\
\label{invariance6}\LL(V_{ij})(V, U) &=& \sum_{pq} (VU)_{ip}U^{*}_{qj}\LL(V_{pq}).
\eeqna

\elem

In the following we give the description of all the diffusion operators of the matrix processes regarding the polar decomposition of $\bfm$.

\bthm\label{th.polar}
\benum
\item{
The action of the  generator on  $\bfH$ and  $\bfN$  is given by
\beqna
\label{aeq1} \Gamma(H_{ij}, H_{kl}) &=& 2(\delta_{jk}H_{il} + \delta_{il}H_{jk}),\\ 
\label{aeq2} \LL(H_{ij}) &=& 4d\delta_{ij},\\
\label{aeq3}\Gamma(N_{ij}, N_{kl}) &=& \sum_{rs} 2\frac{x^2_{r} + x^2_{s}}{(x_{r} + x_{s})^2}U_{ir}\bar{U}_{js}U_{ks}\bar{U}_{lr},\\ 
\label{aeq4}  \LL(N_{ij}) &=&  4\sum_{rs}\frac{x_s}{(x_r + x_s)^2}U_{ir}\bar{U}_{jr}.
\eeqna
Let $\bfD = \{x_1, \cdots, x_{d}\}$ with all $\{x_{i}\}$ different and non negative, then the action of the generator is  given by
\beqna
\label{aeq5}\Gamma(x_{i}, x_{j}) &=& \delta_{ij},\\
\label{aeq6}\LL(x_{i}) &=& \frac{1}{x_{i}} + 4x_{i}\sum_{s \neq i}\frac{1}{x^{2}_{i} - x^{2}_{s}}.
\eeqna
}
\item{
The action of the  generator on  $\bfU$ iss given by
\beqna
\label{aeq7}\Gamma(U_{ij}, U_{kl}) &=&  -r_{lj}U_{il}U_{kj}, \\ 
\label{aeq8}\Gamma(U_{ij}, \bar{U}_{kl}) &=& \delta_{lj}\sum_{s}r_{js}U_{is}\bar{U}_{ks},\\
\label{aeq9}\LL(U_{ij}) &=& -U_{ij}\sum_{s \neq j}r_{js}, \\ 
\label{aeq10}\LL(\bar{U}_{ij}) &=& -\bar{U}_{ij}\sum_{s \neq j}r_{js},
\eeqna
where  $r_{ij} =  2\frac{x^2_{i} + x^2_{j}}{(x^2_{i} - x^2_{j})^2}$.

If we write $\bfW = \bfV \bfU$ and denote $\omega_{ij} = -r_{ij}1_{i \neq j} - \frac{1}{x^{2}_i}\delta_{ij}$, the action of the  generator on  $\bfW$ is given by
\beqna
\label{aeq11}\Gamma(W_{ij}, W_{kl}) &=& \omega_{jl}W_{il}W_{kj},\\  
\label{aeq12}\Gamma(W_{ij}, \bar{W}_{kl}) &=& -\delta_{jl}\sum_{s}\omega_{js}W_{is}\bar{W}_{ks},\\
\label{aeq13}\LL(W_{ij}) &=& \sum_{s \neq j} \omega_{sj}W_{ij}, \\
\label{aeq14}\LL(\bar{W}_{ij}) &=& \sum_{s \neq j} \omega_{sj}\bar{W}_{ij}.
\eeqna

Moreover, 
\beqna
\label{aeq15}\Gamma(W_{ij}, x_{k}) = 0, \\
\label{aeq16}\Gamma(U_{ij},x_{k}) = 0.
\eeqna
}
\eenum
\ethm

\bpf 
Since $\bfH = \bfm^{*}\bfm$, direct computations yield \eqref{aeq1}, \eqref{aeq2}. 
We first compute the action of the  generator on  $\bfD$. Let $P(X) = \det(X\Id - \bfH)$ be the characteristic polynomial and denote the eigenvalues of $\bfH$ by $\{X_{i}\}^{d}_{i=1}$, where $X_{i} = x^{2}_{i}$, then from \eqref{aeq1}, \eqref{aeq2}, we derive
\beqnas
\Gamma(\log P(X), \log P(Y)) &=& \frac{4}{Y-X}(\frac{XP'(X)}{P(X)}-\frac{YP'(Y)}{P(Y)}),\\
\LL(\log P) &=& -4X\frac{P'^2}{P^2}.
\eeqnas
Compare it with the formulas in terms of the eigenvalues, we obtain
\beqnas
\Gamma(X_i, X_j) = 4X_i\delta_{ij}, \ \ \ \LL(X_{i}) = 4(1+ 2X_i\sum_{i\neq j}\frac{1}{X_i-X_j}),
\eeqnas
which lead to \eqref{aeq5}, \eqref{aeq6}.

We now compute the action of the  generator on $\bfU$. Notice that we only need to compute $\Gamma$ and $\LL$ at $\bfU = \Id$, since they are left invariant.
First we have
$$\Gamma(H_{ij}, \log P(X))= \sum^{d}_{r =1} \frac{1}{X_r-X}\Gamma(H_{ij}, X_r)= \sum^{d}_{p, q=1} (H-X)^{-1}_{qp}\Gamma(H_{ij}, H_{pq}).$$
Setting $V_{ijp,k}= \Gamma(U_{ip}\bar U_{jp}, X_k)$, we may derive from the above formula that
$$4 U\frac{D}{D-X} U^*+ \sum_{kp} \frac{X_p}{X_k-X} V_{ijp,k} = 4\frac{H}{H-X} ,$$
which leads to
$$ \sum_{kp} \frac{X_p}{X_k-X} V_{ijp,k}=0.$$
Therefore,
\beqnas
\Gamma(H_{ij}, X_{k}) &=& 4U_{ik}X_{k}\bar{U}_{jk},\\
\sum_{p}X_{p}V_{ijp, k} &=& 0.
\eeqnas
On the other hand, since $\bfU$ is unitary, we have $\sum_{p} \Gamma(U_{ip}\bar{U}_{jp}, \cdot) = 0$, which leads to 
\beqna
\label{aeq17}\Gamma(U_{ij}, \cdot)(\bfU = \Id) = -\Gamma(\bar{U}_{ji}, \cdot)(\bfU = \Id).
\eeqna
Then taking $\sum_{p}X_{p}V_{ijp, k} = 0$ at $U = \Id$, we derive
$$(X_{i} - X_{j})\Gamma(U_{ij}, X_{k}) = 0,$$
which indicates that $\Gamma(U_{ij}, X_{k}) = 0$, for any $i \neq j$ and $k$.

Computing $\Gamma(H_{ij}, H_{kl})$ at $\bfU = \Id$ leads to
\beqnas
2\delta_{il}\delta_{kj}(X_i+X_j) &=& 4 \delta_{i=j=k=l}X_i + X_iX_k\Gamma(\bar U_{ji},\bar U_{lk})\\
& &+ X_iX_l\Gamma(\bar U_{ji}, U_{kl})+ X_jX_l\Gamma(U_{ij}, U_{kl})+ X_jX_k\Gamma(U_{ij}, \bar U_{lk}).
\eeqnas

Then by \eqref{aeq17}, we obtain
$$
2\delta_{il}\delta_{kj}(X_i+X_j)=4 \delta_{i=j=k=l}X_i +(X_i-X_j)(X_k-X_l)\Gamma(U_{ij},U_{kl}),$$
from which we may deduce that for  $i\neq j$ or $k \neq l$
\beqna\label{zeq21}
\Gamma(U_{ij},U_{kl})(\Id) =  -r_{ij}\delta_{il}\delta_{kj},
\eeqna
where $r_{ij} = 2\frac{(X_i+X_j)}{(X_i-X_j)^2}$.

Observe that  $\Gamma(U_{ii}, U_{ii})= r_{ii}$ does not play any role in the computations. This is due to the fact that  $(\bfU,\bfD) \rightarrow \bfU \bfD \bfU^*$ is not a local homeomorphism, and the choice of $U$ is only unique up to a phase. In fact, let $\bfP$ be the unitary matrix of the form $\diag\{e^{i\phi_1}, ..., e^{i\phi_d}\}$ for $0 \leq \phi_1 \leq \cdots \leq \phi_d \leq 2\pi$, then if $\bfU$ is the unitary part of the spectral decomposition of $\bfH$, so is $\bfU \bfP$.  Therefore, we may choose $\bfU$ such that its diagonal elements are all real. Then by \eqref{aeq17}, at $\bfU = \Id$ we have
$$
\Gamma(U_{ii}, \cdot) = 0,
$$
hence $r_{ii} = 0$.


As for $\LL(U_{ij})$, by the fact that $\bfU$ is unitary we obtain at $\bfU = \Id$,
$$
\LL(\bar{U}_{ji}) + \LL(U_{ij}) + 2\sum_{r}\Gamma(U_{ir}, \bar{U}_{jr})= 0.
$$
Since $H_{ij} = \sum_{r}U_{ir}X_{r}\overline{U}_{jr}$, we have at $\bfU = \Id$
$$
4d\delta_{ij} = \LL(H_{ij}) = X_{j}\LL(U_{ij}) + X_{i}\LL(\overline{U}_{ji}) + \LL(X_{i})\delta_{ij} + 2\sum_{r}r_{ir}X_{r}\delta_{ij},
$$
which leads to
\beqnas
(X_{j} - X_{i})\LL(U_{ij}) 
&=& 0.
\eeqnas
Since $U_{ii}$ is real, we have
\beqnas
\LL(U_{ij}) = \LL(\bar{U}_{ji}) =  -\sum_{r \neq i}r_{ir}\delta_{ij}.
\eeqnas

We now compute the diffusion on $\bfN$ and $\bfV$. 
It is not difficult to check that
$$
\Gamma(N_{ij},N_{kl})= \delta_{ijkl}+ \I_{i\neq j} 2\delta_{il}\delta_{jk} \frac{x_i^2+ x_j^2}{(x_i+x_j)^2}=2 \delta_{il}\delta_{kj}\frac{x_i^2+x_j^2}{(x_i+x_j)^2}.
$$

Let $\bfH' = \bfm \bfm^{*}$, then
\beqnas
\Gamma(H'_{ij}, H'_{kl}) &=& 2\delta_{il}H'_{kj} + 2\delta_{jk}H'_{il},\\
\LL(H'_{ij}) &=& 4d\delta_{ij},
\eeqnas
which is exactly the same as $\bfH$. Since $\bfH' = \bfW \bfD^{2} \bfW^{*}$, following the same procedure we may compute $\Gamma$ and $\LL$ of $\bfW$, which is in the same form as $\bfU$. Therefore at $\bfW = \Id$, for $i \neq j$, $k \neq l$ we have
$$
\Gamma(W_{ij}, W_{kl}) = -2\frac{(x^{2}_i+x^{2}_j)}{(x^{2}_i-x^{2}_j)^2}\delta_{il}\delta_{kj}.
$$
and also $\Gamma(W_{ij}, x_k) = 0$.

Since $m_{ij}= \sum^{d}_{p=1}W_{ip}\lambda_p\bar{U}_{jp}$ and
$\Gamma(m_{ij}, m_{kl})=0$,  $\Gamma(m_{ij}, \bar m_{kl})= 2\delta_{ik}\delta_{jl}$, we obtain
\beqnas
0 &=& \Gamma(W_{ij}, W_{kl})x_jx_l + \Gamma(x_i, x_k)\delta_{ij}\delta_{kl} + \Gamma(U_{ij}, U_{kl})x_ix_k\\
& & - \Gamma(W_{ij}, U_{kl})x_jx_k - \Gamma(W_{kl}, U_{ij})x_ix_l,\\
2\delta_{ik}\delta_{jl} &=&  -\Gamma(W_{ij}, W_{lk})x_jx_l + \Gamma(x_i, x_k)\delta_{ij}\delta_{kl} - \Gamma(U_{ij}, U_{lk})x_ix_k\\
& & + \Gamma(W_{ij}, U_{lk})x_jx_k + \Gamma(W_{lk}, U_{ij})x_ix_l,
\eeqnas
which leads to
$$
\Gamma(U_{ij}, W_{kl}) = - \frac{4x_ix_j}{(x^{2}_i - x^{2}_j)^{2}}\delta_{il}\delta_{kj},
$$
and when $i = j$ or $k = l$,
$$
\Gamma(W_{ij}, W_{kl}) = - \frac{1}{x^{2}_i}\delta_{ijkl}.
$$

We may compute $\Gamma$ and $\LL$ of $\bfV$ from $\bfW  = \bfV \bfU$ and $\bfU$. Notice that at $\bfV = \bfU = \Id$, 
\beqnas
\Gamma(W_{ij}, W_{kl}) &=& \Gamma(V_{ij}, V_{kl}) + \Gamma(U_{ij}, U_{kl}) + \Gamma(V_{ij}, U_{kl}) + \Gamma(V_{kl}, U_{ij}),\\
\Gamma(U_{ij}, W_{kl}) &=& \Gamma(U_{ij}, V_{kl}) + \Gamma(U_{ij}, U_{kl}),
\eeqnas
which leads to
\beqnas
\Gamma(U_{ij}, V_{kl}) &=& \frac{2}{(x_i + x_j)^{2}}\delta_{il}\delta_{kj},\\
\Gamma(V_{ij}, V_{kl}) &=& - \frac{4}{(x_i + x_j)^{2}}\delta_{il}\delta_{kj}.
\eeqnas

Also by $\bfN = \bfU \bfD \bfU^{*}$, we have
\beqnas
\Gamma(V_{ij}, N_{kl}) &=& 2\frac{(x_{i} - x_{j})}{(x_{i} + x_{j})^2}\delta_{il}\delta_{jk} := S_{ij}\delta_{il}\delta_{jk},
\eeqnas
and
\beqnas
\Gamma(\bar{V}_{ij}, N_{kl}) &=& 2\frac{(x_{i} - x_{j})}{(x_{i} + x_{j})^2}\delta_{ik}\delta_{jl},
\eeqnas
which is due to $\Gamma(\bar{V}_{ij}, \cdot) = -\Gamma(V_{ji}, \cdot)$ at $\bfV = \Id$. Then at $\bfU = \Id$, 
\beqnas
\LL(N_{ij}) &=& (\LL(x_{i}) - 2\sum_{s}r_{is}x_{s} + 2\sum_{s}x_{i}r_{is})\delta_{ij}\\
&=& 4\sum^{d}_{s=1}\frac{x_s}{(x_i + x_s)^2}\delta_{ij}.
\eeqnas
Now we may compute $\LL(V_{ij})$ and  $\LL(\bar{V}_{ij})$ from $\LL(m_{ij})$. At $\bfV= \bfU=\Id$, we have
\beqnas
\LL(m_{ij}) &=&  \sum_{s}V_{is}\LL(N_{sj})+\LL(V_{is})N_{sj}+ 2\Gamma(V_{is}, N_{sj})\\
&=& 4\sum_{s}\frac{x_s}{(x_i + x_s)^2}\delta_{ij} + \LL(V_{ij})x_{j} + 4\sum_{s}2\frac{(x_{i} - x_{s})}{(x_{i} + x_{s})^2}\delta_{ij},
\eeqnas
such that
\beqnas
\LL(V_{ij}) = \LL(\bar{V}_{ij}) =  -4\sum_{s} \frac{1}{(x_i + x_s)^2}\delta_{ij}.
\eeqnas

Now we have $\Gamma$ and $\LL$ for all the elements in the polar decomposition of the complex matrix $\bfm$ at $\bfV = \bfU = \Id$.
Then by Lemma \ref{invariance}, we are able to obtain all the formulas at $(\bfV, \bfU)$, which ends the proof.

\epf

\section*{Acknowledgements}
This work was done when I was a PhD student at Universit\'e Paul Sabatier. I would like to thank my supervisor, Prof. D. Bakry, for his precious suggestions, helpful discussions and all the encouragements.

\end{document}